\documentclass[amstex,12pt,thmsa,sw20lart]{article}
\usepackage{amsfonts, amsmath,amsthm,amssymb,color}

\input amssym.def
\input amssym.tex

\begin{document}

\date{}
\title{{\bf Pointwise multipliers of \\
Calder\'{o}n-Lozanovski{\u \i} spaces}}

\author{Pawe{\l } Kolwicz\thanks{Research partially supported by the State Committee 
for Scientific Research, Poland, Grant N N201 362236}\,, \ Karol Le\'{s}nik$^\ast$ \,{\small and} 
Lech Maligranda}

\date{}
\maketitle

\vspace{-7mm}

\begin{abstract}

\noindent {\footnotesize Several results concerning multipliers of symmetric Banach function 
spaces are presented firstly. Then the results on multipliers of Calder\'{o}n-Lozanovski{\u \i} 
spaces are proved. We investigate assumptions on a Banach ideal space $E$ and three 
Young functions $\varphi_1, \varphi_2$ and $\varphi $, generating the corresponding 
Calder\'{o}n-Lozanovski{\u \i} spaces $E_{\varphi_1}, E_{\varphi_2}, E_{\varphi}$ so that the 
space of multipliers $M(E_{\varphi_1}, E_{\varphi})$ of all measurable $x$ such that 
$x\,y \in E_{\varphi}$ for any $y \in E_{\varphi_1}$ can be identified with $E_{\varphi_2}$.
Sufficient conditions generalize earlier results by Ando, O'Neil, Zabre{\u \i}ko-Ruticki{\u \i}, 
Maligranda-Persson and Maligranda-Nakai. There are also necessary conditions 
on functions for the embedding $M(E_{\varphi_1}, E_{\varphi}) \subset E_{\varphi_2}$ 
to be true, which already in the case when $E = L^1$, that is, for Orlicz spaces 
$M(L^{\varphi_1}, L^{\varphi}) \subset L^{\varphi_2}$ give a solution of a
problem raised in the book \cite{Ma89}. Some properties of a generalized
complementary operation on Young functions, defined by Ando, are
investigated in order to show how to construct the function $\varphi_2$ such that  
$M(E_{\varphi_1}, E_{\varphi}) = E_{\varphi_2}$. There are also several examples 
of independent interest.}
\end{abstract}

\renewcommand{\thefootnote}{\fnsymbol{footnote}}


\vspace{-5mm}
\begin{center}
\textbf{1. Introduction and preliminaries}
\end{center}

Pointwise multiplication and the space of pointwise multipliers between
Orlicz spaces as well as between some other Banach ideal spaces were
investigated by several authors. Here we try to prove such theorems for the
Calder\'{o}n-Lozanovski\u{\i} spaces $E_{\varphi}$ generated by the Banach
ideal space $E$ and the Young function $\varphi$, which are generalizations
of Orlicz spaces, Orlicz-Lorentz spaces and contain the $p$-convexification $
E^{(p)} (1 \leq p < \infty)$ of $E$. The spaces $E_{\varphi}$ were introduced by
Calder\'on \cite[p. 122]{Cal64} and Lozanovski{\u \i} \cite{Lo65} (see also
Lozanovski{\u \i} \cite{Lo73}). Geometry of the spaces $E_{\varphi}$ was
intensively investigated during the last 20 years (see, for example, \cite
{KL10} and the references given there) and we should also mention here that
they are, in fact, special cases of general Calder\'{o}n-Lozanovski{\u \i}
spaces $\rho(E, F)$ for $F = L^{\infty}$, being important in the
interpolation theory (cf. \cite{KPS}, \cite{Ma89}).

Let $(\Omega, \Sigma, \mu)$ be a complete $\sigma$-finite measure space and $
L^0=L^0(\Omega)$ be the space of all classes of $\mu$-measurable real-valued
functions defined on $\Omega$. A Banach space $E=\left( E,\| \cdot
\|_{E}\right) $ is said to be a \textit{Banach ideal space} on $\Omega$ if $
E $ is a linear subspace of $L^0(\Omega)$ and satisfies the so-called ideal
property, which means that if $y\in E, x \in L^{0}$ and $|x(t)| \leq| y(t)|$
for $\mu $-almost all $t \in \Omega$, then $x\in E$ and $\| x\|_{E} \leq \|
y\| _{E}$.

A Banach ideal space $E$ on $\Omega$ is \textit{saturated} if every $A \in
\Sigma$ with $\mu(A) > 0$ has a subset $B \in \Sigma$ of finite positive
measure for which $\chi_B \in E$. For any such space $E$ it is possible to
construct a set $\Omega_E \in \Sigma$ such that: (i) every element of $E$
vanishes $\mu$-a.e. on $\Omega \setminus \Omega_E$ and (ii) every measurable 
$A \subset \Omega_E$ with $\mu(A) > 0$ has a measurable subset $B$ of finite
positive measure with $\chi_B \in E$. Furthermore, $\Omega_E$ is the union
of an expanding sequence of sets $\{A_k\}$ such that $\mu(A_k) < \infty$ and $
\chi_{A_k} \in E$ for each $k \in \mathbb{N}$. A set $\Omega_E$ is called the 
{\it support of E} and denoted by ${\it supp} E$. Note that we should say here
``a support" rather than ``the support" since in general there will be other
sets ${\tilde \Omega_E}$ which can also satisfy (i) and (ii). However, they
coincide $\mu$-a.e with $\Omega_E$, that is, $\mu({\tilde \Omega_E}
\setminus \Omega_E) = \mu(\Omega_E \setminus {\tilde \Omega_E}) = 0$. It is
also clear that any Banach ideal space $E$ can always be naturally
identified with a saturated Banach ideal space on a possibly smaller
measure space $\Omega_E$. In such space $E$ there exists an element $x_0$
which is strictly positive $\mu$-a.e. on $\Omega_E$, for example, $x_0 =
\sum_{k=1}^{\infty} \chi_{A_k}/(2^k \| \chi_{A_k} \|_E)$. In particular, for
a Banach ideal space $E$ we have ${\it supp} E = \Omega$ if and only if $E$ has a 
\textit{weak unit}, i.e., a function $x$ in $E$ which is positive $\mu $-a.e.
on $\Omega$ (see \cite{KA77} and \cite{Ma89}).

A point $x \in E$ is said to have \textit{order continuous norm} if for any
sequence $(x_n)$ in $E$ such that $0 \leq x_n \leq |x|$ and $x_n \rightarrow
0 ~ \mu$-a.e. on $\Omega$ we have $\| x_n \|_E \rightarrow 0$. By $E_a$ we
denote the subspace of all order continuous elements of $E$. It is known
that $x \in E_a$ if and only if $\| x \chi_{A_n} \|_E \downarrow 0$ for any
sequence $\{A_n\}$ satisfying $A_n \searrow \emptyset$ (that is $A_n \supset
A_{n+1}$ and $\mu(\bigcap_{n=1}^{\infty} A_n) = 0$). A Banach ideal space $E$
is called \textit{order continuous} if every element of $E$ has order
continuous norm, that is, $E = E_a$.

We say that $E$ has the \textit{Fatou property} if $0 \leq x_n \uparrow
x \in L^0$ with $x_n \in E$ and $\sup_{n \in \mathbb{N}} \|x_n\|_E < \infty$
imply that $x \in E$ and $\|x_n\|_E \uparrow \|x\|_E$.

If we consider the space $E$ over a non-atomic measure $\mu$ with 
${\it supp} E = \Omega$, then we say that $E$ is a \textit{Banach function space}. 
If we replace the measure space $\left( \Omega,\Sigma ,\mu \right) $ by the counting 
measure space $\left( \mathbb{N}, 2^{\mathbb{N}}, m\right) ,$ then we say that 
$E$ is a \textit{Banach sequence space} (denoted by $e$). In the last case
the symbol $e_{k}=\left( 0, \ldots,0,1,0, \ldots \right)$ stands for the
k-th unit vector.

The \textit{weighted Banach function space} $E(w)$, where $w$ is a measurable
positive function (weight) on $\Omega$, is defined by the norm $\| x
\|_{E(w)} = \| x w \|_E$.

More information about Banach function spaces and Banach sequence spaces can
be found, for example, in \cite{BS88}, \cite{KA77}, \cite{KPS} and \cite {LT79}.

\medskip Let $E$ and $F$ be ideal Banach spaces in $L^0(\Omega)$ with their
norms $\|\cdot\|_E$ and $\|\cdot\|_F$, respectively. The space of pointwise
multipliers $M(E, F)$ is defined by 
$$
M(E, F) = \{x \in L^0(\Omega): x y \in F ~\text{for all} ~y \in E \}
$$
with the usual operator norm. This space is important, for example, in
investigation of superposition operators and in factorization theorems. Some
properties of superposition operators may be expressed by means of
multiplicator spaces (cf. \cite{AZ87}, \cite{AZ90}). They are also appearing in
factorization theorems. Lozanovski{\u \i} proved that every function $x\in
L^{1}$ can be factorized by $y\in E$ and $z\in E^{\prime }$ in such a way
that $x=yz$ and $\| y\| _{E}\| z\| _{E^{\prime }}\leq
(1+\varepsilon ) \| x\| _{L^{1}}$, where $\varepsilon > 0$ is an arbitrary 
number (cf. \cite{Lo69}). This theorem can be generalized to the form 
$F = E \cdot M(E,F)$ under some assumptions on the spaces (see \cite{Re81}, 
\cite{Sc10}). In the case of sequence spaces (not necessarily ideal) the spaces 
$M(E,F)$ were investigated in \cite{AS76} and used for description of different 
spaces of analytic functions on the disk by sequence multipliers of Taylor 
coefficients. More details about the space $M(E,F) $ we put in the next section.

\medskip 
In this paper we give improvements of the results on multipliers
known for Orlicz spaces $L^{\varphi}$ to the more general situation of 
Calder\'{o}n-Lozanovski{\u \i} spaces $E_{\varphi}$. We need to recall some
necessary definitions about Orlicz and Calder\'{o}n-Lozanovski{\u \i} spaces.

\medskip A function $\varphi :[0,\infty )\rightarrow \lbrack 0,\infty ]$ is
called a {\it Young function} (or \textit{Orlicz function} if it is
finite-valued) if $\varphi $ is convex, non-decreasing with $\varphi (0)=0$;
we assume also that $\varphi $ is neither identically zero nor identically
infinity on $(0,\infty )$ and 
$\lim_{u \rightarrow b_{\varphi}^-} \varphi(u) = \varphi (b_{\varphi})$ if $b_{\varphi} < \infty$, 
where $b_{\varphi} = \sup\{u > 0: \varphi(u) < \infty \}$.

Note that from the convexity of $\varphi $ and the equality $\varphi (0)=0$ it
follows that $\lim_{u\rightarrow 0+}\varphi (u)\newline = \varphi (0)=0$. Furthermore,
from the convexity and $\varphi \not\equiv 0$ we obtain that $
\lim_{u\rightarrow \infty }\varphi (u)=\infty $.

If we denote $a_{\varphi }=\sup \{u\geq 0:\varphi (u)=0\}$, then 
$0 \leq a_{\varphi} \leq b_{\varphi} \leq \infty$ and $
a_{\varphi} < \infty, ~ b_{\varphi} > 0$, since a Young function is neither
identically zero nor identically infinity on $(0,\infty)$. Moreover, $
a_{\varphi} = 0$ if $\varphi$ is $0$ only at $0$ and $b_{\varphi} = \infty$
if $\varphi(u) < \infty$ for $u \in [0, \infty)$. If $\varphi$ takes only
two values $0$ and $\infty$, then $0 < a_{\varphi} = b_{\varphi} < \infty$.
The function $\varphi$ is continuous and nondecreasing on $[0,b_{\varphi})$
and is strictly increasing on $[a_{\varphi}, b_{\varphi})$.

For a given Banach ideal space $E$ on $\Omega$ and a Young function $\varphi$
we define on $L^{0}(\Omega)$ a convex semimodular $I_{\varphi }$ by 
\begin{equation*}
I_{\varphi }(x): = \left\{ 
\begin{array}{cc}
\| \varphi \circ |x| \| _{E} & \text{if }\varphi \circ |x | \in E, \text{ }
\\ 
\infty & \text{otherwise,}
\end{array}
\right.
\end{equation*}
where $(\varphi \circ | x |) (t) =\varphi ( |x(t)| ), t\in \Omega$. By the
Calder\'{o}n-Lozanovski{\u \i} space $E_{\varphi }$ we mean 
\begin{equation*}
E_{\varphi} = \{x \in L^{0}:I_{\varphi }(cx)<\infty ~{\rm  for ~ some } ~ c = c(x) >0\},
\end{equation*}
which is a Banach ideal space on $\Omega$ with the so-called \textit{Luxemburg-Nakano 
norm} defined by 
\begin{equation*}
\| x\| _{E_\varphi }=\inf \left\{ \lambda >0: I_{\varphi}\left( x/\lambda
\right) \leq 1\right\} .
\end{equation*}
If $E = L^1$ ($E = l^1$), then $E_{\varphi}$ is the Orlicz function
(sequence) space $L^{\varphi}$ ($l^{\varphi}$) equipped with the
Luxemburg-Nakano norm (cf. \cite{KR61}, \cite{Ma89}). If $E$ is a Lorentz 
function (sequence) space $\Lambda_w$ ($\lambda_w$), then $E_{\varphi}$ 
is the corresponding Orlicz-Lorentz function (sequence) space $\Lambda_{\varphi, w}$ 
($ \lambda_{\varphi, w}$), equipped with the Luxemburg-Nakano norm. On the
other hand, if $\varphi(u) = u^p, 1 \leq p < \infty$, then $E_{\varphi}$ is
the $p$-convexification $E^{(p)}$ of $E$ with the norm $\|x\|_{E^{(p)}} = \|
|x|^p\|_E^{1/p}$. If $\varphi(u) = 0$ for $0 \leq u \leq 1$ and $\varphi(u)
= \infty$ for $u > 1$, then $E_{\varphi} = L^{\infty}$ with equality of
the norms.
If ${\it supp} E=\Omega $, then ${\it supp} E_{\varphi }=\Omega$, that is, $E_{\varphi }$ 
has a weak unit.

For two ideal Banach spaces $E$ and $F$ on $\Omega $ the symbol $E\overset{C}{
\hookrightarrow }F$ means that the embedding $E\subset F$ is continuous with
the norm which is not bigger than C, i.e., $\Vert x\Vert _{F}\leq C\Vert
x\Vert _{E}$ for all $x\in E$. In the case when the embedding $E\overset{C}{
\hookrightarrow }F$ holds with some (unknown) constant $C>0$ we simply write 
$E\hookrightarrow F$. Moreover, $E = F$ (and $E\equiv F$) means that the spaces
 are the same and the norms are equivalent (equal).

The paper is organized as follows: In Section 1 some necessary definitions
and notation are collected, including the Calder\'{o}n-Lozanovski\u{\i}
spaces $E_{\varphi }$. In Section 2 the space of pointwise multipliers $M(E,F)$
is defined and some general results are presented. In Theorem 1,
some important results in the case of symmetric spaces $E,F$ on $[0,1]$ and $
[0,\infty )$ are proved. It is important to mention here that for symmetric
spaces on $[0,1]$ we have that $M(E,F)\neq \{0\}$ if and only if we have
the imbedding $E\hookrightarrow F$. Also the fundamental function of $M(E,F)$ is
described in terms of fundamental functions $f_{E}$ and $f_{F}$. Better
results appeared in two cases, when either as $E$ we have the smallest symmetric
space (the Lorentz space $\Lambda _{f_{E}}$) or when $E$ is the largest
Marcinkiewicz space $M_{\phi _{1}}$ and $F$ is the smallest Lorentz
space $\Lambda _{\phi }$. Section 3 contains information about the Young
function and its relations with its inverse. Then three relations between
three Young functions are defined and some results proved for two of these
relations (relations for large and small arguments).

Section 4 investigates the embedding $E_{\varphi _{2}}\hookrightarrow
M\left( E_{\varphi _{1}}, E_{\varphi }\right) $. In Theorem 2 there are
sufficient conditions on the Young functions $\varphi _{1},\varphi _{2},\varphi $
and on the Banach ideal space $E$ for such an inclusion. In Theorem 3 and 4
there are necessary conditions on functions under some additional
assumptions on the space $E$. In the special case when $E=L^{1}$ and the
corresponding spaces are Orlicz spaces, then these theorems where proved 
already by Ando \cite{An60} and O'Neil \cite{ON65}.

Section 5 deals with a more difficult reverse embedding $M\left( E_{\varphi
_1}, E_{\varphi }\right) \hookrightarrow E_{\varphi _2} $. The special
case $M(L^p, L^1) \hookrightarrow L^{p^{\prime}}$ is the famous Landau
resonance theorem, which was extended to the case $M(L^{\varphi}, L^1)
\hookrightarrow L^{\varphi^*}$ by several authors (see, for example, \cite
{MW91}, where there are results for Orlicz space $L^{\varphi}$ being even a 
quasi-Banach space). The first results on the embedding $M(L^{\varphi_1},
L^{\varphi}) \hookrightarrow L^{\varphi_2}$, that is, for Orlicz spaces
generated by Orlicz functions on non-atomic measure space, were given by
Zabre{\u \i}ko-Ruticki{\u \i} \cite{ZR67} and Maligranda-Persson \cite{MP89}.
Using a recent result of Maligranda and Nakai \cite{MN} for Orlicz spaces
on general $\sigma$-finite measure spaces and for arbitrary Young functions
we were able to adopt this proof to the situation of Calder\'{o}n-Lozanovski\u{\i} 
spaces $E_{\varphi }$ (Theorem 5). Theorem 6 is interesting here
since under certain monotonicity assumption it was possible to get also
a necessary condition on Young functions for the embedding $M\left(
E_{\varphi _1}, E_{\varphi }\right) \hookrightarrow E_{\varphi _2}$. This
result, for the special case of Orlicz spaces in which case $E = L^1$, gives
an answer to the problem posed in Maligranda's book \cite[Problem 4, p. 77]
{Ma89} under additional assumption of monotonicity of ratio of the fundamental
functions $f_{E_{\varphi _1}}$ and $f_{E_{\varphi}}$.

In Section 6 we have collected, as corollaries from some results in Sections 4 and 5,
the necessary and sufficient conditions on functions so that the equality 
$M\left( E_{\varphi _1}, E_{\varphi }\right) = E_{\varphi _2} $ holds
provided $E$ is a Banach ideal space with the Fatou property and 
${\it supp} E = \Omega$.
2
Section 7 contains construction of a new function from two Young functions
defined probably for the first time by Ando \cite{An60}. This is a
complementary function to $\varphi_1$ with respect to $\varphi$ given by the
formula 
$$
\left( \varphi \ominus \varphi _{1}\right) \left( u\right) =\sup_{v
> 0}\left [ \varphi \left( uv\right) -\varphi _{1}\left( v\right) \right].
$$
The result on this construction gave possibility to improve Theorem 6 having
another monotonicity condition (Theorem 7). Finally, in Example 8 we
show that this last monotonicity condition cannot be dropped. This Example 8
presents construction of an Orlicz function $\psi$ such that the non-separable
Orlicz space $L^{\psi}[0, 1]$ is a proper subspace of $L^2[0, 1]$ and $
M(L^{\psi}[0, 1], L^2[0, 1]) = L^{\infty}[0, 1]$. Moreover, the space $
L^{\psi}[0, 1]$ is also not $L^2[0, 1]$-perfect, that is, $(L^{\psi})^{L^2
L^2}: = M(M(L^{\psi}, L^2), L^2) \neq L^{\psi}$. This is because $(L^{\psi})^{L^2 L^2} = 
(L^{\infty})^{L^2} = L^2 \neq L^{\psi}$.

\vspace{3mm}
\begin{center}
\textbf{2. On the space of pointwise multipliers $M(E, F)$}
\end{center}

Let $E$ and $F$ be ideal Banach spaces in $L^{0}(\Omega )$ with their norms $
\| \cdot \| _{E}$ and $\|  \cdot \| _{F}$, respectively. The
space of pointwise multipliers $M(E,F)$ is defined by 
\begin{equation}
M(E, F)=\{x\in L^{0}(\Omega ):xy\in F~\text{for all}~y\in E\}  \label{multip}
\end{equation}
and the functional on it 
\vspace{-2mm}
\begin{equation}
\| x\| _{M(E, F)}=\sup \{\| xy\| _{F},~y\in E,\| y\|_{E} \leq 1\}  \label{multip-norm}
\end{equation}
defines a complete semi-norm. It is a norm and $M(E,F)$ is an ideal Banach
space if and only if ${\it supp} E = \Omega$, that is, $E$ has a {\it weak unit}, 
i. e., $x_{0}\in E$ such that $x_{0}>0~\mu $-a.e. on $\Omega $ (in
particular, $E\neq \{0\}$). In the case when $F=L^{1}$ we have $
M(E,L^{1})=E^{\prime }$, where $E^{\prime }$ is the classical associated
space to $E$ or the K\"{o}the dual space of $E$, and which is a Banach
function space provided ${\it supp} E = \Omega $. Moreover, ${\it supp} 
E^{\prime }\subset {\it supp} E$ and they are equal if $
\Vert \cdot \Vert _{E}$ has the Fatou null property (if $x_{n}\uparrow x$
and $\Vert x_{n}\Vert _{E}=0$ for all $n\in \mathbb{N}$, then $\Vert x\Vert
_{E}=0$). Always $E\overset{1}{\hookrightarrow }E^{\prime \prime }$ and $
E\equiv E^{\prime \prime }$ if and only if $E$ has the Fatou property. Note that 
$M(E,F)$ can be $\{0\}$ and it can be that ${\it supp} M(E, F)$ is smaller
than  ${\it supp} E \cap {\it supp} F$ (cf. Example 1(c) below).

The notation $E^{\prime}$ for the associated space to $E$ is the reason why
sometimes the space $M(E, F)$ is denoted as $E^F$. Banach ideal spaces for
which $E \equiv E^{\prime \prime}$ are sometimes called perfect spaces and
therefore the Banach ideal space $E$ is called \textit{$F$-perfect} if $E
\equiv E^{FF}$. For example, $L^{\infty}$ and $F$ with ${\it supp} F = \Omega$
are $F$-perfect. Also $E^F$ is $F$-perfect provided ${\it supp} F = {\it supp} E^F =
\Omega$ and $E$ is $L^1$-perfect if and only if $E$ has the Fatou property.

\medskip
General properties and several calculated concrete examples can be found in 
\cite{AZ90}, \cite{MP89} \cite{Ru79} (see also \cite{AS76}, \cite{CDS08}, 
\cite{Cr72}, \cite{CN03}, \cite{Ma89}, \cite{MN}, \cite{Na95}, \cite{Sc10} and \cite{Za66}). 
Let us collect some of these properties and examples:

(i) If $E_0 \overset{C}{\hookrightarrow }E_1$, then $M(E_1, F) \overset{C}{
\hookrightarrow }M(E_0, F)$.

(ii) If $F_0 \overset{C}{\hookrightarrow }F_1$, then $M(E, F_0) \overset{C}{
\hookrightarrow }M(E, F_1)$.

(iii) $E \overset{1}{\hookrightarrow }E^{FF}$ and this embedding follows
from the H\"older-Rogers inequality 

\hspace{7mm} of the form

\hspace{10mm} $\| x y\|_F \leq \| x\|_E \cdot \sup_{\| z \|_E \leq 1} \| y
z\|_F = \| x\|_E \cdot \| y\|_{M(E, F)} $

\hspace{6mm} for any $x \in E$ and $y \in M(E, F)$.

(iv) $E^F$ is $F$-perfect, that is, $E^F \equiv E^{F F F}$.

(v) The embedding $L^{\infty} \overset{C}{\hookrightarrow }M(E, F)$ holds if
and only if $E \overset{C}{\hookrightarrow }F$.

(vi) If ${\it supp} E = \Omega$, then $M(E, E) \equiv L^{\infty}$.

(vii) $M(E, F) \overset{1}{\hookrightarrow }M(F^{\prime}, E^{\prime}) \equiv
M(E^{\prime \prime}, F^{\prime \prime})$. If $F$ has the Fatou property, then

\hspace{7mm} $M(E, F) \equiv M(F^{\prime}, E^{\prime})$.

(viii) $M(E, F) \overset{1}{\hookrightarrow }M(F^G, E^G) \equiv M(E^{G G}, F^{G G})$. 
If $F$ is $G$-perfect, then

\hspace{8mm} $M(E, F) \equiv M(E^{G G}, F^{G G})$.

(ix) For $1 < p < \infty$ we have $M(E^{(p)}, F^{(p)}) \equiv M(E, F)^{(p)}$.

(x) We have equality $\| x \|_{M(E(w_1), F(w_2))} \equiv \| \frac{w_2}{w_1}
x \|_{M(E, F)}$ for $x \in M(E(w_1), F(w_2))$.

\hspace{6mm} In particular, if ${\it supp} E = \Omega $, then $M(E(w_1), E(w_2))
\equiv L^{\infty}(w_2/w_1)$.

(xi) If $F$ has the Fatou property, then $M(E, F)$ also has this property.

\vspace{3mm} 
\textbf{Example 1}. (a) If $1 \leq q < p < \infty, 1/r = 1/q - 1/p, E$ has
the Fatou property and supp\,E = $\Omega$, then $M(E^{(p)}, E^{(q)}) \equiv
E^{(r)}$. In particular, $M(E^{(p)}, E) \equiv E^{(p^{\prime})}$, $
M(L^p(\mu), L^q(\mu)) \equiv L^r(\mu)$ and $M(L^p_{w_1}, L^q_{w_2}) \equiv
L^r_{w_2/w_1}$ for $1 \leq q \leq p \leq \infty$.

(b) Let $1 \leq p < q < \infty$. If the measure $\mu$ is non-atomic, then $
M(L^p(\mu), L^q(\mu)) = \{ 0\}$. Moreover, $M(l^p, l^q) \equiv l^{\infty}$.

(c) Let $\Omega = [0, 2]$ with the Lebesgue measure m. If $E = L^1[0, 1]
\oplus L^2[1, 2]$ with $\|x\|_E = \| x\|_{L^1[0,1]} + \| x\|_{L^2[1, 2]}$
and $F = L^2[0, 2]$, then $M(E, F) = L^{\infty}[1, 2]$. 

(d) Let $\Omega = [0, 2]$ with the Lebesgue measure m. If $E = L^1[0, 1]
\oplus L^{\infty}[1, 2]$ with $\|x\|_E = \| x\|_{L^1[0,1]} + \| x\|_{L^{\infty}[1, 2]}$, 
then $E_a = L^1[0, 1], {\it supp} E_a = [0, 1], {\it supp} E = [0, 2]$ 
and for any Young functions $\varphi, \varphi_1$ we have $E_{\varphi} =
L^{\varphi}[0, 1] \oplus L^{\infty}[1, 2]$ and $M(E_{\varphi_1}, E_{\varphi}) 
= M(L^{\varphi_1}[0, 1], L^{\varphi}[0, 1]) \oplus
L^{\infty}[1, 2]$.

\medskip 
We also need some results in the case of symmetric spaces. By
a {\it symmetric function space} (symmetric Banach function space) on $
I $, where $I = [0, 1]$ or $I = [0, \infty)$ with the Lebesgue measure $m$,
we mean a Banach ideal space $E = (E, \| \cdot\|_E)$ with the additional
property that for any two equimeasurable functions $x \sim y, x, y \in
L^0(I) $ (that is, they have the same distribution functions $d_x = d_y$,
where $d_x(\lambda) = m(\{t \in I: |x(t)| > \lambda \}), \lambda \geq 0$)
and $x \in E$ we have that $y \in E$ and $\| x \|_E = \| y \|_E$. In
particular, $\| x\|_E = \| x^\ast\|_E$, where $x^\ast (t)=\mathrm{inf}
\{\lambda > 0\colon\ d_x(\lambda)\leq t\},\ t \geq 0$.

The {\it fundamental function} $f_{E}$ of a symmetric function space $E$
on $I$ is defined by the formula $f_{E}(t)=\Vert \chi _{\lbrack 0,\,t]}\Vert
_{E},t\in I$. It is well-known that each fundamental function is
quasi-concave on $I$, that is, $f_{E}(0)=0,f_{E}(t)$ is positive, non-decreasing
and $f_{E}(t)/t$ is non-increasing for $t\in (0,m(I))$ or, equivalently, 
$f_{E}(t)\leq \max (1,t/s)f_{E}(s)$ for all $s,t\in (0,m(I))$. Taking ${
\tilde{f}_{E}}(t):= \inf_{s\in (0,m(I))}(1+\frac{t}{s})f_{E}(s)$ we obtain
that the function ${\tilde{f}_{E}}$ is concave and $f_{E}(t)\leq {\tilde{f}_{E}}(t)\leq
2f_{E}(t)$ for all $t\in I$. For any quasi-concave function $\phi $ on $I$
the {\it Marcinkiewicz function space} $M_{\phi }$ is defined by the norm 
\begin{equation*}
\| x \| _{M_{\phi }}=\sup_{t\in I}\phi (t)\,x^{\ast \ast }(t),~~x^{\ast
\ast }(t) = \frac{1}{t}\int_{0}^{t}x^{\ast }(s) ds.
\end{equation*}
This is a symmetric Banach function space on $I$ with the fundamental function 
$f_{M_{\phi }}(t)=\phi (t)$ and $E\overset{1}{\hookrightarrow }M_{f_{E}}$
since 
\begin{equation}
x^{\ast \ast }(t)\leq \frac{1}{t}\, \| x^{\ast } \|
_{E}\Vert \chi _{\lbrack 0,t]}\Vert _{E^{\prime }}=\Vert x\Vert _{E}\frac{1}{
f_{E}(t)}~\mathrm{for~any}~t\in I.  \label{nier1}
\end{equation}
The fundamental function of a symmetric function space $E=(E, \| \cdot \|_{E})$ 
is not necessary concave but we can introduce an equivalent norm on $E$
in such a way that the fundamental function will be concave. In fact, for
the fundamental function $f_{E}$ of $E$ consider the new norm on $E$ defined
by formula 
\begin{equation*}
\Vert x\Vert _{E}^{1}=\max (\Vert x\Vert _{E},\Vert x\Vert _{M_{{\tilde{f}
_{E}}}}),~x\in E.
\end{equation*}
Then $\Vert x\Vert _{E}\leq \Vert x\Vert _{E}^{1}\leq \max (\Vert x\Vert
_{E},2\Vert x\Vert _{M_{f_{E}}})\leq 2\Vert x\Vert _{E}$. Moreover, 
\begin{equation*}
\Vert \chi _{\lbrack 0,t]}\Vert _{E}^{1}=\max (f_{E}(t),{\tilde{f}_{E}}(t))={
\tilde{f}_{E}}(t)
\end{equation*}
and $(E,\Vert \cdot \Vert _{E}^{1})$ is a symmetric Banach function space
with concave fundamental function (cf. Zippin \cite{Zi71}, Lemma 2.1).

For any symmetric function space $E$ with concave fundamental function $f_E$
there is also a smallest symmetric space with the same fundamental
function. This space is the \textit{Lorentz function space} given by the
norm 
\begin{equation*}
\| x \|_{\Lambda_{f_E}} = \int_I x^*(t) d f_E(t) = f_E(0^+) \| x
\|_{L^{\infty}(I)} + \int_I x^*(t) f_E^{\prime}(t) dt.
\end{equation*}
We have then embeddings 
\begin{equation*}
\Lambda_{f_E} \overset{1}{\hookrightarrow }E \overset{1}{\hookrightarrow }
M_{f_E},
\end{equation*}
and all fundamental functions are $f_E$.

Any non-trivial symmetric function space $E$ on $I$ ($E$ is non-trivial if $
E \neq \{0\}$) is intermediate space between the spaces $L^1(I)$ and $
L^{\infty}(I)$. More precisely, 
\begin{equation*}
L^1(I) \cap L^{\infty}(I) \overset{C_1}{\hookrightarrow }E \overset{C_2}{
\hookrightarrow }L^1(I) + L^{\infty}(I),
\end{equation*}
where $C_1 = 2 f_E(1), C_2 = 1/f_E(1)$ and $\|x \|_{L^1 \cap L^{\infty}} =
\max (\|x \|_{L^1}, \|x \|_{L^{\infty}} )$, $\|x \|_{L^1 +L^{\infty}} = \inf
\left\{ \| x_0\|_{L^1} + \|x_1 \|_{L^{\infty}}: x = x_0 +x_1, x_0 \in L^1,
x_1 \in L^{\infty}\right\} = \int_0^1 x^*(s) ds$ (see \cite{KPS}, Theorem
4.1). In particular, supp $E = I$.

A symmetric function space $E$ on $I$ has the \textit{majorant property} if
for all $x \in L^0, y \in E$, the condition $\int_0^t x^*(s)\, ds \leq
\int_0^t y^*(s)\, ds$ for all $t \in I$ implies that $x \in E$ and $\|x \|_E
\leq \| y \|_E$. Every symmetric function space with the Fatou property or
separable symmetric function space have the majorant property. More
information about symmetric spaces on $I = [0, \infty)$ can be found in the
book \cite{KPS}.

\vspace{3mm} 
\textbf{Theorem 1}. \label{symm+fundamental} {\it Let $E$ and $F$ be
non-trivial symmetric function spaces on $I$. 
\vspace{-2mm} 

\begin{itemize}
\item[$(i)$] Then the space of multipliers $M(E, F)$ is a symmetric function space
on $I$. 
\vspace{-1mm} 

Moreover, if the symmetric spaces $E, F$ are on $I =
[0, 1]$, then $M(E, F) \neq \{0\}$ if and only if $E \hookrightarrow F$.
\vspace{-2mm} 

\item[$(ii)$] If $F$ has the majorant property, then $M(E, F)$ has
also majorant property and 
\begin{equation} \label{equation6}
\| x \|_{M(E, F)} = \sup_{\| y\|_E \leq 1} \| x^* y^* \|_F.
\end{equation}
\vspace{-7mm} 

\item[$(iii)$] We have $f _{M(E, F)} (t) \geq \sup_{0 < s\leq t} \frac{f_F(s)}{f_E(s)}$ for all 
$0 < t < m(I)$, and if $f_{F}$ is a concave function with $f_{F}(0^+) = 0$, then
\begin{equation} \label{equation7}
f _{M(E, F)} (t) \leq \int_0^t \frac{f_F^{\prime}(s)}{f_E(s)} \, ds ~{\rm for ~ all} ~ t \in (0, m(I)).
\end{equation}
If, in addition,  $\frac{f_{F}(t)}{f_{E}(t)\, t^{a}}$ is a non-decreasing function on $(0, b)$
for some $a>0$ and $b\in (0, m(I))$, then 
\begin{equation} \label{equation8}
\frac{f_{F}(t) }{f_{E}(t) }\leq f _{M(E, F)} (t) \leq \frac{1}{a} \, \frac{
f_F(t) }{f _{E}(t) } ~~ \mathrm{for ~all} ~ t \in (0, b).
\end{equation}
In the case when $f_F$ is only quasi-concave function, then we should multiply the right 
sides of inequality (\ref{equation7}) and (\ref{equation8}) by constant 2. 
\vspace{-4mm} 

\item[$(iv)$] If $f_{E}$ is a concave function on $I$ with $f_{E}(0^+)
= 0$, then 
\begin{equation*}
f _{M(\Lambda_{f_E}, F)} (t) = \sup_{s \leq t} \frac{f_{F}(s)}{f_{E}(s)} ~~ 
\mathrm{for ~all} ~ t \in I,
\end{equation*}
provided the last supremum is finite.
\vspace{-1mm} 

\item[$(v)$] Let $\phi, \psi$ be concave functions on $I$, $\phi(0^+) = \psi(0^+) = 0$ and 
denote $\phi_1(t) = \frac{t}{\psi(t)}$. Then $M_{\phi_1} \overset{C}{\hookrightarrow }\Lambda_{\phi}$ 
with optimal constant $C > 0$ if and only if $C= \int_I \psi^{\prime}(s) \phi^{\prime}(s) \,ds < \infty$. 
Moreover, $M(M_{\phi_1}, \Lambda_{\phi}) \equiv \Lambda_{\eta}$, where 
$\eta(t) = \int_0^t \psi^{\prime}(s) \phi^{\prime}(s) \,ds$ and $\eta(t) < \infty$ for $t \in I$. 
\end{itemize}
}

Before the proof of Theorem 1 we give the embedding results of independent interest. 

\vspace{3mm} 
\textbf{Proposition 1.} \label{Lem1} {\it Let $E$ and $F$ be non-trivial symmetric function spaces 
on $I$.
\vspace{-3mm}

\begin{itemize}

\item[$(i)$] If either $I = [0, 1]$ or  $I = [0, \infty)$ and $M(E, F) \neq \{0\}$, then $\chi_C \in M(E, F)$ 
for each set $C \subset I$ with 
$m(C) < \infty$.
\vspace{-2mm} 

\item[$(ii)$] If $I = [0, 1]$, then $M(E, F) \neq \{0\}$ if and only if $E \hookrightarrow F$. 
\vspace{-2mm} 

\item[$(iii)$] If $I = [0, \infty)$ and there exists $x \in M(E, F)$ with $x^{\ast}(\infty) = 
\lim_{t \rightarrow \infty} x^{\ast}(t) > 0$, then $E\hookrightarrow F$. 
\vspace{-2mm}

\item[$(iv)$] If $I = [0, \infty)$ and $M(E, F) \neq \{0\}$, 
then $E_{\rm fin} \subset  F$, where $E_{\rm fin}$ is the space of elements from $E$ 
with supports having finite measure.
\vspace{-2mm} 

\end{itemize}
}

\medskip
{\it Proof of Proposition 1.} (i) Let $0\not=x \in M(E, F)$. Then there are a $\delta > 0$ and 
a measurable set $A \subset I$ with $0 < m ( A) < \infty$ such that 
\begin{equation*}
\delta \, \chi _A \leq | x| \, \chi _A.
\end{equation*}
Thus $\chi_A \in M(E, F)$. We will prove that  $\chi _B \in M(E, F)$ 
for each $B \subset I$ with $m( B) = m( A)$. There is a measure preserving transformation 
$\omega: A \rightarrow B$ such that $\omega(A) = B$ (cf. \cite{Ro88}, Theorem 17, p. 410). 
Then $\chi_A = \chi_B(\omega)$ and for any $y \in E$
$$
y(\omega)\, \chi_A = y(\omega) \,  \chi_B (\omega) \sim y \,  \chi_B, 
$$
where $y(\omega)\, \chi_A$ is a function on $I$ defined as $y(\omega(t))$ if $t \in A$ and 
$0$ if $t \notin A$. 
In fact, for any $\lambda > 0$ we have
\begin{eqnarray*}
m(\{t \in I: |y(w(t))\, \chi_A(t)| > \lambda \}) 
&=&
m(\{t \in A: |y(w(t))| > \lambda \}) \\
&=&
m(\omega^{-1}(\{s \in B: |y(s)| > \lambda \})) \\
&=& 
m(\{s \in B: |y(s)| > \lambda \}) \\
&=& 
m(\{s \in I: |y(s) \, \chi_B(s) | > \lambda \}).
\end{eqnarray*}
Since $\chi_A \in M(E, F)$ and $E$ is symmetric it follows that $y(\omega)\, \chi_A \in F$. Consequently, 
by symmetry of $F$ we have $y \,  \chi_B \in F$. If $C \subset I$, then using the fact that measure is 
non-atomic we can write $C$ as a finite sum of disjoint sets $B_k$ such that $m(B_k) = m(A), k = 1, 2, \ldots, n-1$ 
and $m(B_n) \leq m(A)$. Then $\chi_C = \chi_{\cup_{k=1}^n B_k} = \sum_{k=1}^n \chi_{B_k} \in M(E, F)$.

(ii) The sufficiency follows from general properties (i) and (vi) since the embedding $E \hookrightarrow F$ 
gives that $M(E, F) \hookleftarrow M(F, F) = L^{\infty}$. We need to prove necessity. If $M( E,F) \not =\{ 0\} $, 
then by (i) $\chi_{[0, 1]} \in M(E, F)$. Therefore if $x \in E$, then $x = x \, \chi_{[0, 1]} \in F$.
\medskip

(iii) Let $x \in M(E, F)$ be such that $x^{\ast}(\infty) > 0$. Then the set
$$
A = \{t \in I: |x(t)| \geq x^{\ast}(\infty)/2 \}
$$
has infinite measure. Let $y \in E$ be arbitrary. Assume for the moment that there exists a measure 
preserving transformation $\omega: A \rightarrow I$ such that $\omega(A) = I$. Then
$y(\omega)\, \chi_A \sim y$.
Hence, by the symmetry of $E$, we obtain $y(\omega)\, \chi_A \in E$ and so $x \, y(\omega)\, \chi_A \in F$ 
because $x \in M(E, F)$. However, 
$$
\frac{x^{\ast}(\infty)}{2}\, | y( \omega) | \,  \chi_A \leq | x \, y( \omega) | \,  \chi_A \in F,
$$
that is, $ y( \omega)\, \chi_A \in F$ and also $y \in F$ by symmetry of $F$. The only left is to show that there 
exists a measure preserving transformation $\omega: A \rightarrow I$ such that $\omega(A) = I$. Since 
$m(A) = \infty$ and Lebesgue measure is $\sigma$-finite and non-atomic we can write $A$ as a countable 
sum of disjoint subsets $(A_n)_{n=1}^{\infty}$ of $A$, each  of measure $1$. Using the fact that there are 
measure preserving transformations $\omega_n: A_n \rightarrow [n-1, n)$ such that 
$\omega(A_n) = [n-1, n)$ (cf. \cite{Ro88}, Theorem 17, p. 410) we can define mapping 
$\omega: A \rightarrow I$ such that $\omega(A) = I$ by taking $\omega(t): = \omega_n(t)$ if 
$t \in A_n \,( \,n= 1, 2, \ldots)$. Observe that $\omega$ is a measure preserving transformation 
because for any $B \subset I$ we have
\begin{eqnarray*}
m(\omega^{-1}(B))
&=&
m \left [\omega^{-1} \left ( \bigcup_{n=1}^{\infty}(B\cap[n-1, n))\right) \right] =
m \left[ \bigcup_{n=1}^{\infty} \omega_n^{-1} \left( B \cap[n-1, n)\right) \right] \\
&=&
\bigcup_{n=1}^{\infty} m [ \omega_n^{-1} ( B \cap[n-1, n)) ] = \bigcup_{n=1}^{\infty} m ( B \cap[n-1, n)) = m(B).
\end{eqnarray*}

(iv) If $M(E, F) \neq \{0\}$ and $y \in E_{\rm fin}$, then by Proposition 1(i) we obtain $\chi_{{\rm supp}\, y} \in M(E, F)$, 
which implies $y \in F$.

\vspace{3mm} 
{\bf Remark 1}. Note that $E_{fin}$ needs not be complete therefore the inclusion $E_{fin}\subset F$ is 
not continuous, in general. However, for any $d > 0$ there exists a constant $c = c(d) >0$ such that
$E\mid_{A}\overset{c}{\hookrightarrow } F\mid_{A}$ for each $A \subset I$ with $m(A) = d$, that is, 
$\| x \chi_A \|_F \leq c \| x \chi_A \|_E$ for $x \in E$.

Of course, $E\mid _{A}\overset{c}{\hookrightarrow } F\mid _{A}$ for each $A$
with $m ( A) <\infty $, but we need to show that $c$ depends
only on $d$, not on $A$. Suppose, on the contrary, that there is sequence 
$( A_{n}) $ of sets with $m( A_{n}) =d$ and $E\mid
_{A_{n}}\overset{c_n}{\hookrightarrow } F\mid _{A_{n}}$, where $c_n\rightarrow
\infty $ and $c_{n}$ can not be taken smaller. Choose $c$ so that 
$E\mid _{[ 0, d] }\overset{c}{\hookrightarrow } F\mid _{[ 0, d] }.$
Moreover, since $c_{n}$ were optimal, one can find a sequence 
$(x_{n}) \in \prod E\mid _{A_{n}}$ with $\| x_{n}\|_{E}=1$ such that 
$\frac{c_{n}}{2} \| x_{n} \|_{E}\leq \| x_{n}\| _{F}$.
But for each $n \in \mathbb N$ one has $x_{n}^{\ast }\in E\mid _{\left[ 0,d\right] }$ by
symmetry of $E$ and finally 
\begin{equation*}
\infty \leftarrow \frac{c_{n}}{2} \| x_{n}^{\ast } \| _{E\mid
_{[ 0, d] }}=\frac{c_{n}}{2} \| x_{n} \|_{E\mid
_{A_{n}}}\leq \| x_{n} \| _{F\mid _{A_{n}}} 
= \| x_{n}^{\ast } \| _{F\mid _{[ 0, d] }}\leq c \| x_{n}^{\ast } \| _{E\mid _{[ 0, d] }}=c,
\end{equation*}
and this contradiction proves the claim.

\vspace{3mm} 
{\bf Remark 2}. The proof of the embedding in Proposition 1(i) can also be found in \cite[Theorem 1]{Ab91} 
or \cite[Lemma 5.2]{JMST}, where the authors showed more general results from which, in particular, it 
follows that the existence of a nonzero pointwise multiplier necessarily implies that $E\hookrightarrow F$. 
However, our proof is much simpler. Moreover, the inclusion in Proposition 1(iii) is proved 
in \cite[Corollary 9]{Ab91} and in \cite[Lemma 6.1]{CDS08} but the authors used the fact that $M(E, F)$ 
is already a symmetric space which we want to prove here.

\vspace{3mm} 
\textbf{Example 2}. For symmetric spaces on $I = [0, \infty)$ the relation
$M(E, F) \neq \{0\}$ can happen even if we don't have an embedding $E \subset F$. In fact, 
for $E = L^2$ and $F = L^2 \cap L^1$ on $I = [0, \infty)$ we have $E \not \subset F$ but 
$$
M(E, F) = M(L^2, L^2 \cap L^1) \equiv L^{\infty} \cap L^2.
$$

\medskip
{\it Proof of Theorem 1.} (i) Let $I = [0, 1]$. Assume that $x \sim z$ and $0 \neq z \in M(E, F)$. 
By \cite[Lemma 2.1, p. 60]{KPS} (cf. also \cite{As10}, p. 777) for any $\epsilon>0$ there is 
a measure-preserving mapping $\omega: [0,1]\to [0,1]$ such that
$$
 \|x(\omega) - z \|_{L^\infty}\leq \varepsilon.
$$
Moreover, Proposition 1(ii) guarantees that $E\overset{C}{\hookrightarrow }F$.
Thus, for every $y \in E, \|y\|_E \leq 1$, we have 
 \begin{eqnarray*}
 \| xy \|_F 
 &=& 
 \| x(\omega) y(\omega )\|_F \leq
 \| z y(\omega)\|_F + \| [x(\omega)-z] \, y(\omega)\|_F \\
 &\leq&
\| z y(\omega)\|_F + \varepsilon \| y(\omega)\|_F \leq 
\| z y(\omega)\|_F + C \varepsilon \| y(\omega)\|_E \\ 
&\leq& 
\|z\|_{M(E, F)} + C\varepsilon.
\end{eqnarray*}
Taking the supremum over all such $y$, we obtain
$$
\|x\|_{M(E, F)}\leq \|z\|_{M(E, F)} + C \varepsilon,
$$ 
or, since $\varepsilon>0$ is arbitrary, $\|x\|_{M(E, F)}\leq \|z\|_{M(E, F)}$. 
The reverse inequality can be proved similarly 
and the proof is complete in the case when $I = [0, 1]$.

\medskip
Let $I = [0, \infty)$. We divide the proof into two parts. 

\medskip
\noindent
A. Suppose $x\in M(E,F)$ and $x^{\ast }(\infty )> 0$. Set
\begin{equation*}
x_{1}(t): = \max ( |x(t)|, x^{\ast }(\infty )), ~t \in I.
\end{equation*}
Then, for any $y \in E$, we obtain $x \, y \in F$ and $x^{\ast }(\infty ) \, y \in F$ since by Proposition 1(iii) 
we have imbedding $E \subset F$. Thus $x_1 \, y \in F$ for any $y \in E$ and whence $x_1 \in M(E, F)$.
We will prove that 
\begin{equation} \label{equality7}
\| x \|_{M(E, F)} = \| x_1 \|_{M(E, F)}.
\end{equation}
Clearly, it is enough to show that $ \| x \|_{M(E, F)} \geq \| x_1 \|_{M(E, F)}$. Let $\varepsilon >0$ 
be arbitrary. We can find $y\in E$ with $\| y\|_E \leq 1$ such that
\begin{equation*}
(1 - \varepsilon) \| x_{1} \|_{M(E, F)}  \leq \| x_1 \, y \|_F.
\end{equation*}
Consider two sets
\begin{equation*}
A: = \{ t \in I: (1 - \varepsilon ) \, x^{\ast }(\infty ) \leq |x(t)| \leq (1 + \varepsilon) \, x^{\ast }(\infty )\},
\end{equation*}
and
\begin{equation*}
B: = \{ t \in I: | x(t) | > (1 + \varepsilon ) \, x^{\ast }(\infty )\} .
\end{equation*}
Then
\begin{equation*}
x_1\, \chi _B = | x | \,\chi_B, \, m( A) =\infty, ~{\rm and} ~ x_1 \,\chi_A \geq | x | \, \chi _A \geq ( 1-\varepsilon ) \, x_1 \, \chi_A.
\end{equation*}
Similarly as in the proof of the case (iii) of Proposition 1 we can find measure preserving transformation
$\omega _0: A \rightarrow I \backslash B$ such that $\omega_0(A) =  I \backslash B$. Define 
$\omega: A \cup B \rightarrow I$  by
\begin{equation*}
\omega( t) :=\left\{ 
\begin{array}{cc}
\omega _{0} ( t),  & \text{if } t \in A, \\ 
t, & \text{if } t \in B.
\end{array}
\right. 
\end{equation*}
Then $\omega $ is measure preserving transformation and $\omega( A \cup B ) = I$. Moreover, 
$y( \omega) \, \chi _{A\cup B} \sim y$ 
and
\begin{eqnarray*}
| x \,y( \omega)| \, \chi _{A\cup B} 
&=& 
|x \, y( \omega)| \, \chi_A + |x \, y( \omega)| \, \chi _B \\
&\geq& 
(1 - \varepsilon) \, x_1 \, | y( \omega)| \, \chi_A + x_1 \, | y( \omega) | \, \chi _B \\
&\geq& 
(1 -\varepsilon) \, x_1 \, | y( \omega) | \, \chi_{A \cup B}.
\end{eqnarray*}
On the other hand, since $y( \omega _0) \, \chi _A \sim y \, \chi_{I\backslash B}$ it follows that
\begin{eqnarray*}
d_{x_{1} y}
&=& 
d_{x_1 y\chi _B} + d_{x_1 y \chi _{I\backslash B}} \leq
d_{x_1 y\chi _B} + d_{( 1+\varepsilon) x^{\ast }(\infty )\, y\chi_{I\backslash B}}\\
&=&
d_{x_1 y \chi _B} + d_{( 1+\varepsilon ) x^{\ast }(\infty)\, y( \omega _0) \chi _A} \\
&\leq& 
d_{( 1+\varepsilon) x_1 y( \omega) \chi_B} + d_{( 1+\varepsilon ) x_1 y( \omega _0) \chi_A}\\
&=& 
d_{( 1+\varepsilon) x_1 y( \omega) \chi_{A\cup B}}.
\end{eqnarray*}
Thus, $( x_1 \, y)^{\ast } (t) \leq [ ( 1+\varepsilon) \, x_1 \, y( \omega) \, \chi _{A\cup B} ]^{\ast }(t)$ and
\begin{eqnarray*}
( 1 - \varepsilon ) \| x_1 \|_{M(E, F)} 
&\leq& 
\| x_1 \, y \|_F \leq ( 1 + \varepsilon) \| x_1 \, y( \omega)\, \chi _{A\cup B} \| _F \\
&\leq&
\frac{1 + \varepsilon}{1 - \varepsilon} \,  \| x \, y( \omega)\, \chi _{A\cup B} \| _F \leq
\frac{1 + \varepsilon}{1 - \varepsilon} \,  \| x \| _{M(E, F)},
\end{eqnarray*}
which means that 
$$ 
\| x_1 \|_{M(E, F)} \leq \frac{1 + \varepsilon}{(1 - \varepsilon)^2} \,  \| x \| _{M(E, F)},
$$ 
or, since $\varepsilon > 0$ is arbitrary, $ \| x_1 \|_{M(E, F)} \leq  \| x \| _{M(E, F)}$.

Let $z \sim x$. Then $z_1 \sim x_1$, where $x_1 \in M(E, F)$ and $z_1(t) = \max(|z(t)|, z^*(\infty)), t \in I$. We may follow the proof of (i), applying Lemma 2.1, p. 60 in \cite{KPS}, to conclude that $z_1 \in M(E, F)$ and $\| x_1\|_{M(E, F)} = \| z_1\|_{M(E, F)}$. Clearly, $z_1 \geq |z|$ and so $z \in M(E, F)$. Using then equality (\ref{equality7})  we obtain
$$
\| x \|_{M(E, F)} = \| x_1 \|_{M(E, F)} = \| z_1 \|_{M(E, F)} = \| z \|_{M(E, F)},
$$ 
and symmetry of $M(E, F)$ is proved.
\medskip

B. Assume $x^{\ast }( \infty) =0$ and $0 \neq x \in M(E, F)$. Take any $z \sim x.$ 
Then, by \cite[Lemma 2.1, p. 60]{KPS}, for any $\epsilon>0$ there is a measure-preserving 
mapping $\omega: I \to I$ such that
\begin{equation*}
 \|x - z(\omega) \|_{L^1 \cap L^\infty} < \varepsilon,
\end{equation*}
and, for any $y \in E, \| y \|_E \leq 1$, we have 
 \begin{eqnarray*}
 \| z y \|_F 
 &=& 
 \| z(\omega) y(\omega )\|_F \leq
 \| x y(\omega)\|_F + \| [z(\omega)-x] \, y(\omega)\|_F \\
 &\leq&
\|x\|_{M(E, F)} + \| [z(\omega)-x] \, y(\omega)\|_F.
\end{eqnarray*}
Since $\|y (\omega)\|_{E} \leq 1$, then, by using Lemma 1 proved below, we 
can find a decomposition $y(\omega) = u + v$ such that 
$u \in E, m({\rm supp}\, u) \leq 1$ and $v \in E \cap L^{\infty}$ with $\| v\|_{L^{\infty}} \leq \frac{1}{f_E(1)}$. 
Therefore, applying Proposition 1(iv) on inclusion $E_{\rm fin} \subset  F$, we obtain 
\begin{eqnarray*}
\| z y \|_F 
&\leq& 
\|x\|_{M(E, F)} + \| [z(\omega)-x] \, [u + v] \|_F\\
&\leq&
\|x\|_{M(E, F)} + \| [z(\omega)-x] \, u \|_F +  \| [z(\omega)-x] \, v \|_F\\
&\leq&
\|x\|_{M(E, F)} + \| z(\omega)-x \|_{L^1 \cap L^\infty} \, \|u \|_F +  \| [z(\omega)-x] \, v \|_F\\
&\leq&
\|x\|_{M(E, F)} +\varepsilon \, \|u \|_F + \|v \|_{L^{\infty}} \| z(\omega)-x \|_F\\
&\leq&
\|x\|_{M(E, F)} +\varepsilon \, \|u \|_F + \frac{1}{f_E(1)} \| z(\omega)-x \|_{L^1\cap L^{\infty}} \cdot 2 f_F(1)\\
&=&
\|x\|_{M(E, F)} +\varepsilon \, \|u \|_F + 2 \varepsilon \frac{f_F(1)}{f_E(1)}.
\end{eqnarray*}
Using Remark 1 we have $\| u\|_F \leq c \| u \|_E$ since $m({\rm supp} u) \leq 1$ and, hence,
$$
\| z y \|_F \leq \|x\|_{M(E, F)} +\varepsilon \, c \|u \|_E + 2 \varepsilon \frac{f_F(1)}{f_E(1)}
\leq  \|x\|_{M(E, F)} +\varepsilon \, c + 2 \varepsilon \frac{f_F(1)}{f_E(1)}.
$$
Taking the supremum over all $y \in E, \| y\|_E \leq 1$, we obtain
$$
\| z \|_{M(E, F)} \leq  \|x\|_{M(E, F)} +\varepsilon \, c + 2 \varepsilon \frac{f_F(1)}{f_E(1)},
$$
and since $\varepsilon > 0$ is arbitrary $\| z \|_{M(E, F)} \leq  \|x\|_{M(E, F)}$. The reversed inequality can 
be proved similarly turning the roles of both functions and the proof is complete.

\vspace{3mm} 
\textbf{Lemma 1.} \label{Sublem1} {\it Let $E$ be a symmetric space on $[0, \infty)$. If $y \in E$ and 
$\| y \|_{E} \leq 1$, then we can find a decomposition $y = u + v$ such that $u \in E, \| u \|_E \leq 1$ 
with $m({\rm supp}\, u) \leq 1$ and $v \in E\cap L^{\infty}$ with $\| v\|_{L^{\infty}} \leq 1/f_E(1)$.}

\medskip
{\it Proof of Lemma 1.} If either $y^{\ast}(1) > y^{\ast}(\infty) \geq 0$ or $y^{\ast}(1) = y^{\ast}(\infty)> 0$, then 
we can take $y = u + v$, where $u = y \chi_A$ and $v = y \chi_{I\setminus A}$ with $A = \{s> 0: |y(s)| > y^{\ast}(1) \}$. 
In fact, $m(A) \leq 1$ and 
$$
\| v \|_{L^{\infty}} \leq y^{\ast}(1) \leq \int_0^1 y^{\ast}(t) dt \leq \|y \|_{E} \| \chi_{[0, 1]}\|_{E^{\prime}} \leq 1/f_E(1).
$$
If $y^{\ast} (1) = y^{\ast }( \infty) = 0$ and $A {\rm supp}\,y$, then $m(A) = t$ for some $0 < t \leq 1$. Then 
$y = y \chi_A + 0$ is such a decomposition.

\medskip
{\it Proof of Theorem 1.} (ii) Let $x \in M(E, F), y \in E$. By Theorem 1(i) and symmetry of $E$ we have 
$x^{\ast} \in M(E, F), y^{\ast} \in E$ and thus $x^{\ast} y^{\ast} \in F$. Moreover,  
\begin{eqnarray*}
\sup_{\| y\|_E \leq 1} \| x^* y^*\|_F &=& \sup_{\| y^*\|_E \leq 1} \| x^*
y^*\|_F \\
&\leq& \sup_{\| z\|_E \leq 1} \| x^* z \|_F = \| x^*\|_{M(E, F)}.
\end{eqnarray*}
On the other hand, since a symmetric space $F$ has the majorant property and 
\begin{equation*} 
\int_0^t (x y)^*(s) ds = \sup_{m(A) = t} \int_A |x(s) y(s)| ds \leq \int_0^t
x^*(s) y^*(s) ds ~~ \mathrm{for ~all} ~~ t \in I,
\end{equation*}
it follows that $\| x y \|_F \leq \| x^* y^* \|_F$ and, hence,
\begin{equation*}
\| x \|_{M(E, F)} = \sup_{\| y\|_E \leq 1} \| x y \|_F \leq \sup_{\| y\|_E
\leq 1} \| x^* y^*\|_F \leq \| x^*\|_{M(E, F)} =  \| x\|_{M(E, F)},
\end{equation*}
where the last equality follows from (i), that is, from symmetry of $M(E, F)$
and equality (\ref{equation6}) is proved.

We show that if $F$ has the majorant property, then $M(E, F)$ has it as well. Let $x \in M(E, F), z \in L^0$ 
and 
$$
\int_0^t z^*(s)\, ds \leq \int_0^t x^*(s)\, ds
$$ 
for all $t \in I$. By Hardy lemma 
\begin{equation} \label{8}
\int_0^t z^*(s) \,y^*(s) \, ds \leq \int_0^t x^*(s) \,y^*(s)\, ds
\end{equation}
for all $t \in I$ and each $ y \in E$. Since, by (i), $M(E, F)$ is symmetric it follows that $x^* \in M(E, F)$. 
Then $x^* y^* \in F$ for each $y \in E$. The majorant property of $F$ and inequality (\ref{8}) give $z^* y^* \in F$. 
Following analogously as above and applying majorant property of $F$ we obtain $z y \in F$. Thus, $z \in M(E, F)$ 
and by (\ref{8}) $\| z^* y^* \|_F \leq \| x^* y^* \|_F$. Taking the supremum over all $y \in E$ with $\| y \|_E \leq 1$ 
we have $\| z \|_{M(E, F)} \leq \| x \|_{M(E, F)}$, which shows the majorant property of $M(E, F)$.

\medskip
(iii) For any $t \in I$ we have 
\begin{equation*}
f_{M(E, F)}(t) = \sup_{\| y\|_E \leq 1} \| y \, \chi_{[0, t]} \|_F \geq
\sup_{s \in I} \| \frac{\chi_{[0, s]} }{f_E(s)} \, \chi_{[0, t]} \|_F =
\sup_{s \leq t} \frac{f_F(s)}{f_E(s)}.
\end{equation*}
On the other hand, if $f_F$ is a concave function and $f_F(0^+) = 0$, then
by the general properties (i) and (ii) we have $M(M_{f_E}, \Lambda_{f_F}) \overset{1}{
\hookrightarrow }M(E, F)$ and, hence,
\begin{eqnarray*}
f_{M(E, F)}(t) 
&\leq& 
f_{M(M_{f_E}, \Lambda_{f_F})}(t) = \sup_{\| y
\|_{M_{f_E} \leq 1}} \| y \chi_{[0, t]} \|_{\Lambda_{f_F}} \\
& \leq&
 \sup_{y^* \leq 1/f_E} \int_0^t y^*(s) d f_F(s) \leq \int_0^t \frac{
f_F^{\prime}(s)}{f_E(s)} ds.
\end{eqnarray*}
If, in addition, we have monotonicity assumption on $\frac{f_{F}(t)}{f_{E}(t)\, t^{a}}$, 
then, by using the fact that $f_F^{\prime}(s) \leq f_F(s)/s$ for almost all $s \in I$,  
we obtain for $t \in (0, b)$
\begin{eqnarray*}
 \int_0^t \frac{f_F^{\prime}(s)}{f_E(s)} ds 
&\leq&
\int_0^t \frac{f_F(s)}{s \, f_E(s)} ds = \int_0^t \frac{f_F(s)}{
f_E(s)\, s^a} s^{a-1} ds \\
&\leq& \frac{f_F(t)}{f_E(t)\, t^a} \int_0^t s^{a-1} ds = \frac{1}{a} \, 
\frac{f_F(t)}{\, f_E(t)}.
\end{eqnarray*}

(iv) The estimation follows from the equality $\|x^*\|_{M(\Lambda_{f_E}, F)} =
\sup_{s \in I} \frac{\| x^* \chi_{[0, s]}\|_F}{f_E(s)}$ proved in \cite[Theorem 3]{MP89} since then 
\begin{equation*}
f_{M(\Lambda_{f_E}, F)} (t) = \sup_{s \in I} \frac{\| \chi_{[0, t]} \chi_{[0,
s]}\|_F}{f_E(s)} = \sup_{s \leq t} \frac{f_F(s)}{f_E(s)}.
\end{equation*}

(v) First we show the equivalence. Suppose $\int_I \psi^{\prime}(s) \phi^{\prime}(s) \,ds = C < \infty$. If 
 $\|x\|_{M_{\phi_1}} \leq 1$, then for all $0 < t < m(I)$
 $$
 \frac{1}{\psi(t)} \int_0^t x^*(s) \, ds = \frac{\phi_1(t)}{t}  \int_0^t x^*(s) \, ds \leq 1,
 $$ 
 whence
 $$
 \int_0^t x^*(s) ds \leq \psi(t) = \int_0^t \psi^{\prime}(s) ds
 $$ 
 for each $t \in I$. Thus, by Hardy lemma (see Proposition 3.6, p. 56 in \cite{KPS}), 
 \begin{equation} \label{9}
 \int_I x^*(s) \phi^{\prime}(s) \,ds \leq  \int_I \psi^{\prime}(s) \phi^{\prime}(s) \,ds = C
 \end{equation}
 and so $\| x \|_{\Lambda_{\phi}} \leq C$, which means that $M_{\phi_1} \overset{C}{\hookrightarrow }\Lambda_{\phi}$. 
 Moreover,
 \begin{equation} \label{10}
\| \psi^{\prime} \|_{M_{\phi_1}} = \sup_{t \in I} \phi_1(t) \left(\psi^{\prime} \right)^{**}(t) = 
\sup_{t \in I} \frac{1}{\psi(t)} \int_0^t \psi^{\prime}(s) \,ds = 1.
 \end{equation}
Note that for $x = x^* = \psi^{\prime}$ we have equality in estimate (\ref{9}), that is, $C$ is optimal.
 
Suppose, conversely, $M_{\phi_1} \overset{C}{\hookrightarrow } \Lambda_{\phi}$ with optimal constant $C$. 
 Then
 \begin{equation} \label{11}
\int_I \psi^{\prime}(s) \phi^{\prime}(s) \,ds \leq c \, \| \psi^{\prime} \|_{M_{\phi_1}} = C.
 \end{equation}
Moreover, if $x = x^*$ and $\| x \|_{M_{\phi_1}} = 1$, then $ 1 = \sup_{t \in I} \frac{\phi_1(t)}{t} \int_0^t x(s)\, ds$, 
whence $\int_0^t x(s)\, ds \leq \psi(t) = \int_0^t \psi^{\prime}(s)\, ds$ for any $t \in I$. By Hardy lemma,  
$\int_0^t x(s) \phi^{\prime}(s)\, ds \leq \int_0^t \psi^{\prime}(s) \phi^{\prime}(s) \, ds $ for each $t \in I$. Consequently
$$
\| x \|_{\Lambda_{\phi}} = \int_I x(s) \phi^{\prime}(s)\, ds \leq  \int_0^t \psi^{\prime}(s) \phi^{\prime}(s) \, ds 
$$ 
for any $x \in M_{\phi_1}$ with $\| x \|_{M_{\phi_1}} = 1$. Thus
$$
\| x \|_{\Lambda_{\phi}} \leq  \int_0^t \psi^{\prime}(s) \phi^{\prime}(s) \, ds \, \| x \|_{ M_{\phi_1}},
$$
for all $x \in M_{\phi_1}$. But $ C \leq \int_I \psi^{\prime}(s) \phi^{\prime}(s) \, ds$, because $C$ is optimal, which 
together with (\ref{11}) gives the equality $C = \int_0^t \psi^{\prime}(s) \phi^{\prime}(s) \, ds$.

The equality of spaces follows from the following facts: if $x^*
\in M(M_{\phi_1}, \Lambda_{\phi})$ and since $\| \psi^{\prime}
\|_{M_{\phi_1}} = 1$ we obtain $x^* \psi^{\prime} \in \Lambda_{\phi}$ or 
\begin{equation*}
\int_I (x^* \psi^{\prime})^*(s) \phi^{\prime}(s) \,ds = \int_I x^*(s)\,
\psi^{\prime}(s) \phi^{\prime}(s) \,ds \leq \| x^*\|_{M(M_{\phi_1}, \Lambda_{\phi})},
\end{equation*}
and, thus, $\| x\|_{\Lambda_{\eta}} \leq \| x^*\|_{M(M_{\phi_1}, \Lambda_{\phi})} $. 
On the other hand, if $x \in \Lambda_{\eta}$ and $\|y \|_{M_{\phi_1}} \leq 1$ is 
arbitrary or, equivalently, $\int_0^t y^*(s)\, ds \leq \psi(t) = \int_0^t \psi^{\prime}(s)\, ds$ 
for all $t \in I$, then, by the Hardy inequality, 
\begin{equation*}
\int_I x^*(s) y^*(s) \phi^{\prime}(s) \,ds \leq \int_I x^*(s) \psi^{\prime}(s)
\phi^{\prime}(s) \,ds =\| x\|_{\Lambda_{\eta}}.
\end{equation*}
Thus $\| x y\|_{\Lambda_{\phi}} \leq \| x\|_{\Lambda_{\eta}}$ for any $\|y
\|_{M_{\phi_1}} \leq 1$ and so $\| x^*\|_{M(M_{\phi_1}, \Lambda_{\phi})}
\leq \| x\|_{\Lambda_{\eta}}$. The proof is complete. \qed

\vspace{3mm} 
{\bf Example 3}. A special case of symmetric spaces for which we can calculate the fundamental function 
of their space of multipliers was given in \cite[Example 4.2]{FMP10} and we give a short proof of this result. 
Let $E \hookrightarrow F$ be two ultrasymmetric spaces on $[0, 1]$ with the same parameter $\tilde {G}$, 
that is, 
$$
\| x\|_E = \| f_E(t) x^{\ast}(t)\|_{\tilde {G}}, ~ \| x\|_F = \| f_F(t) x^{\ast}(t)\|_{\tilde {G}},
$$
where $\tilde {G}$ is a symmetric space on $(0, 1)$ with respect to the measure $dt/t$ (see Pustylnik \cite{Pu03} 
and Astashkin-Maligranda \cite{AM08}). Then
\begin{equation} \label{equality8}
f _{M(E, F)} (t) = \sup_{0 < s\leq t} \frac{f_F(s)}{f_E(s)} ~~{\rm for ~ all} ~ t \in (0, 1].
\end{equation}
In fact, for any $\| x\|_E \leq 1$ we have
\begin{eqnarray*}
\| x \chi_{[0, t]} \|_F 
&=&
\| f_F (x \chi_{[0, t]})^{\ast} \|_{\tilde {G}} \leq \| f_F x^{\ast} \chi_{[0, t]} \|_{\tilde {G}} =
\| \frac{f_F}{f_E} f_E x^{\ast} \chi_{[0, t]} \|_{\tilde {G}} \\
&\leq & 
\sup_{0 < s\leq t} \frac{f_F(s)}{f_E(s)} \, \| f_E x^{\ast} \chi_{[0, t]} \|_{\tilde {G}} \leq
  \sup_{0 < s\leq t} \frac{f_F(s)}{f_E(s)},
\end{eqnarray*}
and the reverse estimate is always true by Theorem 1(iii). Thus we obtain equality (\ref{equality8}).
Another example of spaces with equality (\ref{equality8}) will be given in Example 9.

\vspace{3mm}
\begin{center}
\textbf{3. Some properties of Young functions}
\end{center}

To state and prove our main results we will need to define some subclasses
of Young functions, an inverse of Young function and their properties. We
write $\varphi >0$ when $a_{\varphi }=0$ and $\varphi <\infty $ if $%
b_{\varphi }=\infty$. Define the sets of Young functions $\mathcal{Y}%
^{\left( i\right) },$ for $i=1, 2, 3,$ as 
\begin{eqnarray*}
\mathcal{Y}^{\left( 1\right) } &=&\left\{ \varphi : b_{\varphi }=\infty
\right\} ,  \label{klasy} \\
\mathcal{Y}^{\left( 2\right) } &=&\left\{ \varphi : b_{\varphi }<\infty 
\text{ and }\varphi \left( b_{\varphi }\right) =\infty \right\} ,  \notag \\
\mathcal{Y}^{\left( 3\right) } &=&\left\{ \varphi : b_{\varphi }<\infty 
\text{ and } \varphi \left( b_{\varphi }\right) <\infty \right\} .  \notag
\end{eqnarray*}
For an Young function $\varphi$ we define its right-continuous inverse in a
generalized sense by the formula (cf. O'Neil \cite{ON65}): 
\vspace{-2mm}
\begin{equation}  \label{inverse}
\varphi^{-1}(v) = \inf \{u \geq 0: \varphi(u) > v \} ~ \mathrm{for} ~ v \in
[0, \infty) ~~ \mathrm{and} ~~ \varphi^{-1}(\infty) = \lim_{v \rightarrow
\infty} \varphi^{-1}(v).
\end{equation}
Note that $\left\{ u\geq 0:\varphi \left( u\right) >v\right\} \neq \emptyset$
for each $v\in [0, \infty)$.

We will often use properties of an Young function $\varphi$ and its generalized
inverse $\varphi^{-1}$. Therefore let us collect these properties here.

\vspace{3mm} 
\textbf{Lemma 2.} \label{lem-1} \textit{We have 
\vspace{-3mm} }

\begin{itemize}
\item[$(i)$] \textit{$\varphi (\varphi^{-1}(u)) \leq u$ for all $u \in [0,
\infty)$ and $u \leq \varphi^{-1}(\varphi(u))$ if $\varphi (u) <\infty$. 
\vspace{-2mm} }

\item[$(ii)$] \textit{$\varphi ^{-1}(\varphi (u))=u$ for $a_{\varphi }\leq
u\leq b_{\varphi }$ if $b_{\varphi }<\infty $ and $\varphi \left( b_{\varphi
}\right) <\infty $. \vspace{-2mm} }

\item[$(iii)$] \textit{$\varphi ^{-1}(\varphi (u))=u$ for $a_{\varphi }\leq
u<b_{\varphi }$ if either $b_{\varphi }=\infty $ or $b_{\varphi }<\infty $
and $\varphi \left( b_{\varphi }\right) =\infty $. \vspace{-2mm} }

\item[$(iv)$] \textit{$\varphi ^{-1}(\varphi (u))>u$ for $\ 0\leq
u<a_{\varphi }$. \vspace{-2mm} }

\item[$(v)$] \textit{$\varphi ^{-1}(\varphi (u))<u$ for $u>b_{\varphi }$. 
\vspace{-2mm} }

\item[$(vi)$] \textit{$\varphi (\varphi ^{-1}(u))=u$ if $u\in \lbrack
0,\infty )$ and $\varphi \in \mathcal{Y}^{(1)}\cup \mathcal{Y}^{(2)}$. 
\vspace{-2mm} }

\item[$(vii)$] \textit{$\varphi (\varphi ^{-1}(u))=u$ if $u\in \lbrack
0,u_{0}]$ and $\varphi \in \mathcal{Y}^{(3)}$, where $u_{0}=\inf \left\{
u>0:\varphi ^{-1}(u)=b_{\varphi }\right\} $. \vspace{-2mm} }

\item[$(viii)$] \textit{$\varphi (\varphi ^{-1}(u))<u$ if $u>u_{0}$ and 
$\varphi \in \mathcal{Y}^{(3)}$, where $u_{0}=\inf \left\{ u>0:\varphi
^{-1}(u)=b_{\varphi }\right\} $. }
\end{itemize}
Note that (i) follows from (vi)--(viii) and (ii)--(iv).

\medskip 
We will use also the following notations: the symbol $\varphi
_{1}^{-1}\varphi _{2}^{-1}\prec \varphi ^{-1}$ for all arguments [for large
arguments] (for small arguments) means that there is a constant $C>0$ [there
are constants $C,u_{0}>0$] (there are constants $C,u_{0}>0$) such that the
inequality 
\begin{equation}
C\varphi _{1}^{-1} ( u) \varphi _{2}^{-1} ( u) \leq \varphi ^{-1}( u)
\label{left-ineq}
\end{equation}
holds for all $u>0$ [for all $u\geq u_{0}$] (for all $0 < u\leq u_{0}$),
respectively.
\vspace{1mm}

The symbol $\varphi ^{-1}\prec \varphi _{1}^{-1}\varphi _{2}^{-1}$ for all
arguments [for large arguments] (for small arguments) means that there is a
constant $D>0$ [there are constants $D,u_{0}>0$] (there are constants 
$D, u_{0}>0$) such that the inequality 
\begin{equation}
\varphi ^{-1}( u) \leq D\varphi_{1}^{-1}( u) \varphi_{2}^{-1} ( u)
\label{right-ineq}
\end{equation}
holds for all $u>0$ [for all $u\geq u_{0}$] (for all $0 < u\leq u_{0}$),
respectively.

\vspace{2mm}

The symbol $\varphi _{1}^{-1}\varphi _{2}^{-1}\approx \varphi ^{-1}$ for all
arguments [for large arguments] (for small arguments) means that $\varphi
_{1}^{-1}\varphi _{2}^{-1}\prec \varphi ^{-1}$ and $\varphi ^{-1}\prec
\varphi _{1}^{-1}\varphi _{2}^{-1},$ that is provided there are constants 
$C, D > 0$ [there are constants $C,D,u_{0}>0$] (there are constants 
$C, D, u_0 >0$) such that the inequalities
\begin{equation*}
C\varphi _{1}^{-1} ( u) \varphi _{2}^{-1} ( u) \leq \varphi ^{-1} ( u) \leq
D\varphi _{1}^{-1} ( u) \varphi_{2}^{-1} ( u)
\end{equation*}
hold for all $u>0$ [for all $u\geq u_{0}$] (for all $0 < u\leq u_{0}$),
respectively.

\vspace{3mm} 
\textbf{Lemma 3.} \label{lem-2} {\it $( i) $ If $\varphi _{1}^{-1}\varphi
_{2}^{-1}\approx \varphi ^{-1}$ for large arguments, then 
\vspace{2mm}

$( a) $ for any $0<u_{1}<u_{0}$ there are contants 
$C_{1}\leq C, D_{1}\geq D$ such that 
\begin{equation}  \label{C1D1-large}
C_{1} \varphi _{1}^{-1} ( u) \varphi _{2}^{-1} ( u) \leq \varphi ^{-1} ( u)
\leq D_{1} \varphi _{1}^{-1} ( u) \varphi _{2}^{-1} ( u) ~ for ~ any ~ u\geq
u_{1}.
\end{equation}

\vspace{-2mm}

$( b) $ $b_{\varphi }<\infty $ if and only if $b_{\varphi
_{1}}<\infty $ and $b_{\varphi _{2}}<\infty$. 
\vspace{2mm}

\noindent 
$( ii) $ If $\varphi _{1}^{-1}\varphi
_{2}^{-1}\approx \varphi ^{-1}$ for small arguments, then 

\medskip 
$( a) $ for any $u_{1}>u_{0}$ there are contants 
$C_{1}\leq C,D_{1}\geq D$ such that 
\begin{equation}  \label{C1D1-small}
C_{1} \varphi _{1}^{-1} ( u) \varphi _{2}^{-1} ( u) \leq \varphi ^{-1} ( u)
\leq D_{1} \varphi _{1}^{-1} ( u) \varphi _{2}^{-1} ( u) ~ for ~ any ~ 0 <
u\leq u_{1}.
\end{equation}

\vspace{-2mm}

$( b) $ $a_{\varphi }=0$ if and only if $a_{\varphi
_{1}}=0$ or $a_{\varphi _{2}}=0.$ }

\vspace{3mm}
{\it Proof.} (i). In order to prove $\left( a\right)$ it is
enough to take 
\begin{equation*}
C_{1}=\min \left\{ C,\inf\limits_{u_{1}\leq u\leq u_{0}}\frac{\varphi ^{-1}
( u) }{\varphi _{1}^{-1} ( u) \varphi _{2}^{-1} ( u) }\right\} \text{ and} ~
D_{1}=\max \left\{ D,\sup\limits_{u_{1}\leq u\leq u_{0}}\frac{\varphi ^{-1}(
u) }{ \varphi _{1}^{-1} ( u) \varphi _{2}^{-1} ( u) }\right\} .
\end{equation*}
We prove (ii),(b). Necessity. Suppose $a_{\varphi }=0$ and $a_{\varphi
_{1}}>0$ or $a_{\varphi _{2}}>0.$ Taking $u_{n}\rightarrow 0$ we get that
$\varphi _{1}^{-1} ( u_n) \varphi_{2}^{-1} ( u_n) \rightarrow a_{\varphi
_{1}}a_{\varphi _{2}}>0$ and $\varphi ^{-1}\left( u_{n}\right) \rightarrow
0, $ a contradiction with inequality (\ref{left-ineq}). Sufficiency. If $
a_{\varphi _{1}}=0$ and $a_{\varphi }>0,$ then $\varphi _{1}^{-1} ( u_n)
\varphi_{2}^{-1} ( u_n) \rightarrow 0$ and $\varphi ^{-1} (u_n) \rightarrow
a_{\varphi }$, a contradiction with inequality (\ref{right-ineq}). The case 
$a_{\varphi _{2}}=0$ and $a_{\varphi }>0$ can be proved in an analogous way.

The proofs of (i),(b) and (ii),(a) are similar. \qed

\begin{center}
\textbf{4. On the inclusion $E_{\varphi _{2}}\hookrightarrow M\left(
E_{\varphi _{1}}, E_{\varphi }\right) $}
\end{center}


We will consider the question when the product $x y \in E_{\varphi}$ provided $x \in
E_{\varphi_1}$ and $y \in E_{\varphi_2}$.

\vspace{3mm} 
\textbf{THEOREM 2}. \label{thm-1} {\it Suppose $E$ is a Banach ideal space with the Fatou property and 
$\varphi, \varphi _{1}$ and $\varphi _{2}$ are Young functions. Assume also that at least 
one of the following conditions holds: 
\vspace{-3mm} 

\begin{itemize}
\item[$(i)$] $\varphi _{1}^{-1}\varphi _{2}^{-1}\prec \varphi^{-1}$
for all arguments. 
\vspace{-2mm} 

\item[$(ii)$] $\varphi _{1}^{-1}\varphi _{2}^{-1} \prec \varphi^{-1} 
$ for large arguments and $L^{\infty }\hookrightarrow E$. 
\vspace{-2mm} 

\item[$(iii)$] $\varphi _{1}^{-1}\varphi _{2}^{-1} \prec
\varphi^{-1} $ for small arguments and $E\hookrightarrow L^{\infty }$. 
\end{itemize}
\vspace{-2mm}

Then, for every $x \in E_{\varphi_1}$ and $y \in E_{\varphi_2}$ the product 
$x y \in E_{\varphi}$, which means that 
$E_{\varphi _{2}}\hookrightarrow M\left( E_{\varphi _{1}}, E_{\varphi}\right).$ }

\newpage

\textit{Proof.} (i). It is well-known that inequality (\ref{left-ineq}) implies 
\begin{equation}  \label{sufficiency}
\varphi \left( Cuv\right) \leq \varphi _{1}\left( u\right) +\varphi
_{2}\left( v\right)
\end{equation}
for each $u,v>0$ with $\varphi _{1}\left( u\right) ,\varphi _{2}\left(
v\right) <\infty$ (cf. O'Neil \cite{ON65}). In fact, taking $w=\max \left[
\varphi _{1}\left( u\right) ,\varphi _{2}\left( v\right) \right]$, we obtain
\begin{equation*}
uv\leq \varphi _{1}^{-1} ( \varphi _{1} ( u)) \varphi _{2}^{-1} ( \varphi
_{2} ( v)) \leq \varphi _{1}^{-1} ( w) \varphi _{2}^{-1} ( w) \leq \frac{1}{C
} \varphi ^{-1} ( w) ,
\end{equation*}
so that $\varphi \left( Cuv\right) \leq \varphi ( \varphi^{-1} ( w)) \leq
w\leq \varphi _{1} ( u) +\varphi_{2} (v)$. Then, taking any $x\in
E_{\varphi_1}$ and $y \in E_{\varphi_2}$ with $\| x\|_{\varphi _{1}} = \|
y\| _{\varphi _{2}} =1$, we obtain 
\begin{equation*}
I_{\varphi }\left( \frac{Cxy}{2}\right) \leq \frac{1}{2} \,I_{\varphi
}\left( Cxy\right) \leq \frac{1}{2} \, \left[ I_{\varphi _{1}}\left(
x\right) +I_{\varphi _{2}}\left( y\right) \right] \leq 1.
\end{equation*}
This means that $\| xy \|_{E_{\varphi }} \leq \frac{2}{C} \, \|x \|
_{E_{\varphi _{1}}} \,\| y \| _{E_{\varphi _{2}}}$ for any $x\in E_{\varphi
_{1}}$ and $y\in E_{\varphi _{2}}$, and so $E_{\varphi _{2}} \overset{2/C}{
\hookrightarrow }M\left( E_{\varphi _{1}}, E_{\varphi }\right)$.

\medskip (ii) Set $u_{1}=\frac{1}{\Vert \chi _{\Omega }\Vert _{E}}$. Let $
C_{1}=C_{1}\left( u_{1}\right) $ be the corresponding number from (\ref
{C1D1-large}). Then, analogously as in $(i)$, we conclude
\begin{equation*}
\varphi \left( C_{1}uv\right) \leq \varphi _{1}(u)+\varphi _{2}(v)
\end{equation*}
for each $u,v>0$ with $\varphi _{1}(u),\varphi _{2}(v)<\infty $ and $\max
\left\{ \varphi _{1}(u),\varphi _{2}(v)\right\} \geq u_{1}$. Let $x\in
E_{\varphi _{1}},y\in E_{\varphi _{2}}$ with $\Vert x\Vert _{\varphi
_{1}}=\Vert y\Vert _{\varphi _{2}}=1$ and 
\begin{equation*}
A=\left\{ t\in \Omega :\max \left[ \varphi _{1}(|x(t)|),\varphi _{2}(|y(t)|)
\right] \geq u_{1}\right\} ,~\mathrm{and}~~B=\Omega \backslash A.
\end{equation*}
Then 
\begin{equation*}
I_{\varphi }\left( \frac{C_{1}xy}{3}\chi _{A}\right) \leq \frac{1}{3}\,\left[
I_{\varphi _{1}}(x\chi _{A})+I_{\varphi _{2}}(y\chi _{A})\right] \leq \frac{2
}{3}.
\end{equation*}
Since $I_{\varphi _{1}}(x)\leq 1$ it follows that $\varphi
_{1}(|x(t)|)<\infty $ for $\mu $-a.e. $t\in \Omega $ and, consequently, 
\begin{equation*}
|x(t)|\leq \varphi _{1}^{-1}(\varphi _{1}(|x(t)|))\leq \varphi
_{1}^{-1}(u_{1})~\mathrm{for~each}~t\in B.
\end{equation*}
Analogously, $|y(t)|\leq \varphi _{2}^{-1}(u_{1})$. Then, by (\ref{C1D1-large}), 
we obtain 
\begin{equation*}
I_{\varphi }(C_{1}xy\chi _{B})\leq I_{\varphi }(C_{1}\varphi
_{1}^{-1}(u_{1})\varphi _{2}^{-1}(u_{1})\chi _{B})\leq \varphi (\varphi
^{-1}(u_{1}))\Vert \chi _{\Omega }\Vert _{E}\leq u_{1}\Vert \chi _{\Omega
}\Vert _{E}=1.
\end{equation*}
Finally, 
\begin{equation*}
I_{\varphi }\left( \frac{C_{1}xy}{3}\right) \leq I_{\varphi }\left( \frac{
C_{1}xy}{3}\,\chi _{A}\right) +\frac{1}{3}\,I_{\varphi }\left( C_{1}xy\,\chi
_{B}\right) \leq \frac{2}{3}+\frac{1}{3}=1,
\end{equation*}
and, thus, $\Vert xy\Vert _{E_{\varphi }}\leq \frac{3}{C_{1}}.$ Consequently, 
$\| xy\|_{E_{\varphi }}\leq \frac{3}{C_{1}}\| x\| _{E_{\varphi
_{1}}}\| y\| _{E_{\varphi _{2}}}$ for any $x\in E_{\varphi _{1}}$ and $
y\in E_{\varphi _{2}}$. Thus $E_{\varphi _{2}}\overset{3/C_{1}}{
\hookrightarrow }M\left( E_{\varphi _{1}}, E_{\varphi }\right) $.

$(iii)$ Assume that $E\overset{A}{\hookrightarrow }L^{\infty }$. We then
follow analogously as above case $(i)$ showing $E_{\varphi _{2}} \overset{
2/C_{1}}{\hookrightarrow }M\left( E_{\varphi _{1}}, E_{\varphi }\right)$,
where $C_{1} = C_{1}(A) $ is from (\ref{C1D1-small}) for $u_{1}=A$ and $A$
is such that $\mathrm{ess}\sup \limits_{t \in \Omega} |\varphi ( u(t))| \leq
A$ for any $u\in E_{\varphi }$ with $\|u\|_{E_{\varphi }} \leq 1$. \qed

\medskip The next result shows the necessity of the estimate $
\varphi_{1}^{-1}\varphi _{2}^{-1} \prec \varphi ^{-1}$ in Theorem 2.

\vspace{3mm} 
\textbf{THEOREM 3}. \label{thm-2} {\it Let $E$ be a Banach function
space with the Fatou property and let $\varphi, \varphi _{1}, \varphi _{2}$
be Young functions. Suppose 
\begin{equation}  \label{necessity}
E_{\varphi _{2} }\hookrightarrow M \left( E_{\varphi _{1}}, E_{\varphi}\right).
\end{equation}

\vspace{-6mm}

\begin{itemize}
\item[$(i)$] If $E_{a} \not = \{ 0\}$, then $\varphi_{1}^{-1}\varphi
_{2}^{-1} \prec \varphi ^{-1}$ for large arguments. 
\vspace{-2mm} 

\item[$(ii)$] If ${\it supp} E_{a} = {\it supp} E$ and $L^{\infty}\not
\hookrightarrow E$, then $\varphi _{1}^{-1}\varphi _{2}^{-1}\prec \varphi
^{-1}$ for all arguments. 
\end{itemize}
}

{\it Proof.} $(i)$ Suppose the condition $\varphi _{1}^{-1}\varphi
_{2}^{-1}\prec \varphi ^{-1}$ is not satisfied for large arguments. This
means that there exists a sequence $\left( u_{n}\right) $ with 
$u_{n}\nearrow \infty $ such that, for any $n\in \mathbb{N}$,
\begin{equation}
\varphi _{1}^{-1}(u_{n})\varphi _{2}^{-1}(u_{n})\geq 2^{n}\varphi
^{-1}(u_{n}).  \label{neces-9}
\end{equation}
We want to construct a sequence $\left\{ x_n \right\}\subset E_{\varphi
_2} $ such that $\| x_n\| _{E_{\varphi _2}}\leq 1$ but 
$\| x_n\|_{M\left( E_{\varphi _1}, E_{\varphi }\right) }\rightarrow
\infty $, which is equivalent to the fact that $E_{\varphi
_2}\not\hookrightarrow M\left( E_{\varphi _1}, E_{\varphi }\right) $. 
\newline
First of all, note that for almost all $n\in \mathbb N$ we can find
measurable sets $A_{n}$ satisfying 
\begin{equation}
\Vert u_{n}\chi _{A_n}\Vert _E = 1.  \label{neces-10}
\end{equation}
In fact, if $E_a \not=\{0\}$, then there is a nonzero $0\leq x\in E_a$
and therefore there is also a set $A$ of positive measure such that $\chi
_{A}\in E_a$. Of course, for large enough $n$ one has $\left\|
u_n \chi _{A}\right\| _{E}\geq 1.$ Applying Dobrakov result from \cite
{Do74} we conclude that the submeasure $\omega (B)=\| \chi _B \|_{E}$
for $B\in \Sigma, B\subset A,$ has the Darboux property. Consequently, for
each such $n$ there exists a set $A_n$ satisfying (\ref{neces-10}). Define 
\begin{equation*}
x_n = \varphi _{2}^{-1}(u_{n})\chi _{A_n},~~y_n =\varphi_{1}^{-1}(u_{n})\chi _{A_n}.
\end{equation*}
Then 
\begin{equation*}
I_{\varphi _{1}}(y_n) = \| \varphi _1 (\varphi _{1}^{-1}(u_{n})\chi
_{A_n}) \| _{E}\leq \| u_n \chi _{A_n}\| _{E} = 1
\end{equation*}
and thus $\| y_n \| _{E_{\varphi _1}}\leq 1$. Similarly, we can show
that $\|  x_n \| _{E_{\varphi _2}}\leq 1$. However, for large enough 
$n$, one has by (\ref{neces-9}) 
\begin{eqnarray*}
I_{\varphi }\left( \frac{x_n y_n}{\lambda }\right) 
&=&
\| \varphi\left( \frac{\varphi _{1}^{-1}\left( u_{n}\right) \varphi _{2}^{-1}\left(
u_{n}\right) }{\lambda }\right) \chi _{A_{n}}\| _{E} \\
&\geq &
\| \varphi \left( \frac{2^{n}\varphi ^{-1}\left( u_n \right) }{
\lambda }\right) \chi _{A_{n}}\| _{E}.
\end{eqnarray*}
If $\varphi \in \mathcal{Y}^{\left( 1\right) }\cup \mathcal{Y}^{\left(
2\right) }$ and $\lambda \leq 2^{n}$, then 
\begin{eqnarray*}
\| \varphi \left( \frac{2^{n}\varphi ^{-1}(u_{n})}{\lambda }\right) \chi
_{A_{n}}\| _{E} &\geq &\| \frac{2^{n}}{\lambda }\varphi \left( \varphi
^{-1}\left( u_{n}\right) \right) \chi _{A_{n}}\| _{E} \\
&=&\| \frac{2^n}{\lambda }\,u_n \,\chi _{A_{n}}\| _{E}=\frac{2^{n}}{
\lambda }\geq 1.
\end{eqnarray*}
If $\varphi \in \mathcal{Y}^{\left( 3\right) }$, then for sufficiently large 
$n$ and $\lambda <2^n$ we obtain that
\begin{equation*}
I_{\varphi }\left( \frac{x_n y_n}{\lambda }\right) = \| \varphi \left( 
\frac{2^{n}\,b_{\varphi }}{\lambda }\right) \chi _{A_{n}}\| _{E}=\infty ,
\end{equation*}
which implies $\| x_n y_n\|_{\varphi }\geq 2^{n}$. Finally, 
\begin{equation*}
\| x_{n}\| _{M(E_{\varphi _1}, E_{\varphi })}=\sup_{\| y \|_{E_{\varphi _1}}\leq 1}
\| x_n y \| _{E_{\varphi }}\geq \| x_{n}y_{n}\|_{E_{\varphi }}\geq 2^{n},
\end{equation*}
whereas $\| x_n \| _{E_{\varphi _2}}\leq 1$ and this is the required
sequence.

\medskip 
$(ii)$ Of course, the assumption ${\it supp} E_a = {\it supp} E$ implies that 
$E_{a} \not= \{ 0\}$. Therefore we need only to prove that $
\varphi _{1}^{-1}\, \varphi _{2}^{-1} \prec \varphi^{-1}$ for small
arguments. Note that in this case the proof is almost the same as in $(i)$.
One only has to prove that there are sets like in (\ref{neces-10}). Since
${\it supp} E_{a} = {\it supp} E$ we see that there is $x\in E_{a}$ with $x>0$ a.e.
Define $B_{k} = \{ t\in \Omega: x(t) > \frac{1}{k} \}$. Of course, $B_{k}$ have positive 
measure for sufficently large $k$, $\Omega = \bigcup_{k=1}^{\infty} B_{k}$
and $B_{1} \subset B_{2} \subset B_{3} \subset \ldots$. We have $\|
\chi_{B_{k}} \|_{E} \rightarrow \infty$ because $L^{\infty
}\not\hookrightarrow E$ and $E$ has the Fatou property. Moreover, $
\chi_{B_{k}}\in E_{a}$ for any $n \in \mathbb{N}$. Therefore, for each $
u_{n} $ one can find $k\left( n\right)$ such that $\| u_{n}\chi _{B_{k(n)}}
\|_{E} >1$ and the argument is the same as before. \qed

\medskip The following example explains why the conditions concerning $E_{a}$
in Theorem 3 are reasonable but not necessary.

\vspace{3mm} 
\textbf{Example 4}. Note that Theorem 3 is not true without any additional
assumption on the space $E$. In fact, for $E=L^{\infty }$ and for any
non-trivial Young function $\varphi $ we have that $E_{\varphi }=L^{\infty }$,
which gives 
\begin{equation*}
M(E_{\varphi _{1}}, E_{\varphi })=M(L^{\infty }, L^{\infty })=L^{\infty
}=E_{\varphi _{2}},
\end{equation*}
and no relation between the functions $\varphi _1, \varphi _2,\varphi $ is
necessary. On the other hand, for the weighted space $E=L_{t}^{\infty }[0,1]$
with its norm $\| x\| _{E}=~{\rm ess\, sup}_{t\in \left[ 0,1\right]
}\left\vert t\,x(t)\right\vert $ we have $E_{a}=\left\{ 0\right\} $ but the
condition $\varphi _{1}^{-1}\,\varphi _{2}^{-1}\prec \varphi ^{-1}$ for
large arguments is necessary for inclusion $E_{\varphi _{2}}\hookrightarrow
M\left( E_{\varphi _{1}}, E_{\varphi }\right) $. To see this it is enough to
proceed like in Theorem 3(i) because the function $\eta :[0,1]\rightarrow
\lbrack 0,1]$ given by $\eta (t)=\Vert \chi _{\lbrack 0,t]}\Vert _{E}=t$ is
continuous and therefore has the Darboux property.

\vspace{3mm} 
{\bf Remark 3}. The condition $E_{a}\not=\left\{ 0\right\} $ from Theorem
3(i) may be changed by the following weaker one: there is $a>0$ such that for 
any $0<t<a$ we can find $A\in \Sigma $ with $\| \chi _A \| _E = t$.

\vspace{3mm} 
Now we investigate necessity condition on the Young functions
in the case of Banach sequence space.

\vspace{3mm} 
{\bf THEOREM 4}. \label{thm-3} {\it Let $e$ be a Banach sequence
space with the Fatou property and let $\varphi $, $\varphi _{1}$, $\varphi
_{2}$ be Young functions. Suppose 
\begin{equation}
e_{\varphi _{2}}\hookrightarrow M\left( e_{\varphi _{1}}, e_{\varphi }\right).
\label{neces-11}
\end{equation}

\vspace{-6mm}

\begin{itemize}
\item[$(i)$] If $l^{\infty }\not\hookrightarrow e$ and $\sup_{i \in 
\mathbb{N}}\left\Vert e_{i}\right\Vert _{e}<\infty$, then $
\varphi_{1}^{-1}\varphi _{2}^{-1} \prec \varphi ^{-1}$ for small arguments. 
\vspace{-2mm} 

\item[$(ii)$] If $e\not\hookrightarrow l^{\infty }$ and for each $a > 1$ there is a set $B_a$ with 
$\frac{1}{2} \leq a \| \chi_{B_a} \| \leq 1$, then $\varphi_{1}^{-1}\varphi _{2}^{-1} \prec \varphi ^{-1}$ for 
large arguments. 
\end{itemize}
}

\vspace{2mm}

{\it Proof.} $(i)$ Suppose that the condition $\varphi _{1}^{-1}\varphi
_{2}^{-1}\prec \varphi ^{-1}$ is not satisfied for small arguments. Then
there exists a sequence $(u_{n})$ with $u_{n}\rightarrow 0$ such that for
any $n\in \mathbb{N}$ we have
\begin{equation}
\varphi _{1}^{-1}(u_{n})\varphi _{2}^{-1}(u_{n})\geq 2^{n}\varphi
^{-1}(u_{n}).  \label{neces-12}
\end{equation}
Since $l^{\infty }\not\hookrightarrow e$ and $e$
has the Fatou property, it follows that $\| \chi
_{\left\{ 1,2,...,n\right\} }\| _{e}\rightarrow \infty $ as $n\rightarrow
\infty $.
From the assumption $\sup_{i\in \mathbb{N}}\Vert e_{i}\Vert _{e}<\infty $
one can find $N$ large enough such that $\Vert u_{n}\,e_{i}\Vert _{e}\leq
1/2 $ for each $n>N$ and each $i\in \mathbb{N}$. Furthermore, for each $n>N$
we can find $k_{n}$ satisfying 
\begin{equation*}
\left\| u_{n}\chi _{\left\{ 1,2,...,k_{n}\right\} }\right\| _{e}\leq 1
\text{ and }\left\Vert u_{n}\chi _{\left\{ 1,2,...,k_{n},k_{n}+1\right\}
}\right\Vert _{e}>1.
\end{equation*}
Because $\left\| u_{n}e_{k_{n}+1}\right\| _{e}\leq 1/2$ one has also 
$1/2\leq \| u_{n}\chi _{\left\{ 1,2,...,k_{n}\right\} }\| _{e}$. For 
$n\ $large enough we put $A_{n}=\left\{ 1,2,...,k_{n}\right\} $, and 
\begin{equation*}
x_{n}=\varphi _{2}^{-1}(u_{n})\chi _{A_{n}}, ~y_{n}=\varphi
_{1}^{-1}(u_{n})\chi _{A_{n}}.
\end{equation*}
Then $\| y_{n} \| _{e_{\varphi _{1}}}\leq 1$ and $\|x_{n}\|_{e_{\varphi _{2}}}\leq 1$. 
Moreover, by (\ref{neces-12}), one has 
\begin{eqnarray*}
I_{\varphi }\left( \frac{x_{n}y_{n}}{\lambda }\right) &=&\| \varphi
\left( \frac{\varphi _{1}^{-1}\left( u_{n}\right) \varphi _{2}^{-1}\left(
u_{n}\right) }{\lambda }\right) \chi _{A_{n}}\| _{e} \\
&\geq &\| \varphi \left( \frac{2^{n}\varphi ^{-1}\left( u_{n}\right) }{
\lambda }\right) \chi _{A_{n}}\| _{e}\,.
\end{eqnarray*}
If $\lambda \leq 2^{n-1}$, then, by applying Lemma 2(vi) and (vii), we obtain 
\begin{equation*}
\Vert \varphi \left( \frac{2^{n}\varphi ^{-1}\left( u_{n}\right) }{\lambda }
\right) \chi _{A_{n}}\Vert _{e}\geq \Vert \frac{2^{n}u_{n}}{\lambda }\,\chi
_{A_{n}}\Vert _{e}\geq \frac{2^{n-1}}{\lambda }\geq 1,
\end{equation*}
for sufficiently large $n,$ which implies that $\| x_{n}y_{n}\|_{e_{\varphi }}\geq 2^{n-1}$ 
and, consequently, 
\begin{equation*}
\| x_{n}\| _{M\left( e_{\varphi _{1}},e_{\varphi }\right)
}=\sup_{\| y \| _{e_{\varphi _{1}}}\leq 1}\| x_n y \|_{e_{\varphi }}\geq \| x_{n} y_{n}\| _{e_{\varphi }}\geq 2^{n-1},
\end{equation*}
whereas $\|  x_{n}\Vert _{e_{\varphi _{2}}}\leq 1$. Therefore (\ref
{neces-11}) is not satisfied.

\medskip $(ii)$ We proceed as in $\left( i\right) $. Note that the condition is satisfied in many non-symmetric spaces, 
for example, in $e = l^1(\{\frac{1}{i}\})$ for which $l^{\infty} \not\hookrightarrow e$ or $e = l^1(\{\frac{1}{i^2}\})$ 
in which $l^{\infty} \hookrightarrow e$. \qed

\vspace{3mm} 
\noindent 
Putting Theorems 2 and 3 together we obtain:

\vspace{3mm} 
{\bf Corollary 1}. \label{cor-1} {\it Let $E$ be a Banach function
space with the Fatou property and let $\varphi $, $\varphi _{1}$, $\varphi
_{2}$ be Young functions. 
\vspace{-2mm} 

\begin{itemize}
\item[$(i)$] Suppose $L^{\infty }\hookrightarrow E$ and $E_{a}\not=
\{ 0 \}$. Then $E_{\varphi _{2}}\hookrightarrow M\left( E_{\varphi
_{1}}, E_{\varphi }\right) $ if and only if $\varphi_{1}^{-1}\varphi
_{2}^{-1}\prec \varphi ^{-1}$ for large arguments. 
\vspace{-2mm} 

\item[$(ii)$] Assume $L^{\infty }\not\hookrightarrow E$ and ${\it supp} E_{a} = {\it supp} E$. 
Then $E_{\varphi _{2}}\hookrightarrow M\left(
E_{\varphi _{1}}, E_{\varphi }\right) $ if and only if $\varphi_{1}^{-1}
\varphi _{2}^{-1}\prec \varphi ^{-1}$ for all arguments. 
\end{itemize}
}

Before giving a similar characterization for the sequence case note that the equality
$e = l^{\infty }$ implies then $e_{\varphi _{2}} = e_{\varphi _{1}} = e_{\varphi }= 
M\left(e_{\varphi_{1}}, e_{\varphi }\right) =l^{\infty }$ for any Orlicz functions.
Consequently, looking for a neccesary and sufficient condition for the
inclusion $e_{\varphi _{2}}\hookrightarrow M\left( e_{\varphi
_{1}},e_{\varphi }\right) $ we need to consider the following three cases:
(i) $e\hookrightarrow l^{\infty }$ and $l^{\infty}\not\hookrightarrow e$,
(ii) $l^{\infty }\hookrightarrow e$ and $e\not\hookrightarrow l^{\infty }$,
(iii) $l^{\infty }\not\hookrightarrow e$ and $e\not\hookrightarrow l^{\infty
}.$

\noindent Taking into account Theorems 2 and 4, we then get immediately

\vspace{3mm} 
{\bf Corollary 2}. \label{cor-2} {\it Let $e$ be a Banach sequence
space with the Fatou property and let $\varphi, \varphi _{1}, \varphi _{2}$
be Young functions. 
\vspace{-2mm} 

\begin{itemize}
\item[$(i)$] Suppose $e\hookrightarrow l^{\infty }, l^{\infty
}\not\hookrightarrow e$ and $\sup_{i \in \mathbb{N}}\ \| e_{i} \| _{e}<
\infty$. Then $e_{\varphi _{2}}\hookrightarrow M( e_{\varphi
_{1}}, e_{\varphi}) $ if and only if $\varphi _{1}^{-1}\varphi _{2}^{-1}\prec
\varphi ^{-1}$ for small arguments.
 \vspace{-2mm} 
\end{itemize} 
Suppose additionally that for each $a > 1$ there is a set $B_a$ with $\frac{1}{2} \leq a \| \chi_{B_a} \| \leq 1$.
\begin{itemize}
\item[$(ii)$] Assume that $l^{\infty }\hookrightarrow e$ and 
$e\not\hookrightarrow l^{\infty }$. Then $e_{\varphi _{2}}\hookrightarrow
M\left( e_{\varphi _{1}}, e_{\varphi }\right) $ if and only if $\varphi
_{1}^{-1}\varphi _{2}^{-1}\prec \varphi ^{-1}$ for large arguments. 
\vspace{-2mm} 

\item[$(iii)$] Let $l^{\infty }\not\hookrightarrow e,$ 
$e\not\hookrightarrow l^{\infty }$ and $\sup_{i\in \mathbb{N}} \| e_{i}
\|_{e}<\infty$. Then $e_{\varphi _{2}}\hookrightarrow M\left( e_{\varphi
_{1}}, e_{\varphi }\right) $ if and only if $\varphi
_{1}^{-1}\varphi_{2}^{-1}\prec \varphi ^{-1}$ for all arguments. 
\end{itemize}
}

Note that conditions (\ref{left-ineq}) and (\ref{sufficiency}) are
equivalent. The implication (\ref{left-ineq}) $\Rightarrow $ (\ref
{sufficiency}) was shown in the proof of Theorem 2. We prove that $\varphi
(Cuv)\leq \varphi _{1}(u)+\varphi _{2}(v)$ for all $u,v>0$ implies $\varphi
_{1}^{-1}(w)\varphi _{2}^{-1}(w)\leq \frac{2}{C}\varphi ^{-1}(w)$ for all $
w>0$. In fact, for $w>0$ let $u=\varphi _{1}^{-1}(w)$ and $v=\varphi
_{2}^{-1}(w)$. Then, by assumption and Lemma 2, we have 
\begin{equation*}
\varphi (Cuv)\leq \varphi _{1}(u)+\varphi _{2}(v)=\varphi _{1}(\varphi
_{1}^{-1}(w))+\varphi _{2}(\varphi _{2}^{-1}(w))\leq 2w,
\end{equation*}
and again, by Lemma 2 and the concavity of $\varphi ^{-1}$, we obtain 
\begin{equation*}
Cuv\leq \varphi ^{-1}(\varphi (Cuv))\leq \varphi ^{-1}(2w)\leq 2\varphi
^{-1}(w),
\end{equation*}
which gives $\varphi _{1}^{-1}(w)\varphi _{2}^{-1}(w)\leq \frac{2}{C}\varphi
^{-1}(w)$. Now we consider the respective case for large and small
arguments. Discussing the case for large arguments we prove that the
following conditions are equivalent:

$\left( a\right) $ for any $u_{0}>0$ there is $C_{0}>0$ such that 
$C_{0}\varphi _{1}^{-1}(u)\varphi _{2}^{-1}(u)\leq \varphi ^{-1}(u)$ for all 
$u\geq u_{0}.$

$\left( b\right) $ for any $u_{1}>0$ there is $C_{1}>0$ such that $\varphi
(C_{1}uv)\leq \varphi _{1}(u)+\varphi _{2}(v)$ for all $u,v$ 

\hspace{5mm} with $u_{1}\leq
\max \left\{ \varphi _{1}(u),\varphi _{2}(v)\right\} <\infty .$

\noindent
The implication $\left( a\right) \Rightarrow \left( b\right) $ has been
shown in the proof of Theorem 2(ii). We prove $\left( b\right) \Rightarrow
\left( a\right).$ For any Orlicz function $\varphi \in \mathcal{Y}^{\left(
3\right) }$ denote 
\begin{equation*}
\alpha _{\varphi }=\inf \left\{ u>0:\varphi ^{-1}(u)=b_{\varphi }\right\}
\end{equation*}
and for $\varphi \in \mathcal{Y}^{\left( 1\right) }\cup \mathcal{Y}^{\left(
2\right) }$ we set $\alpha _{\varphi }=\infty .$ Let
\begin{equation*}
u_{1}=\min \left\{ u_{0},\alpha _{\varphi _{1}},\alpha _{\varphi
_{2}}\right\} .
\end{equation*}
Take $w\geq u_{0}$ and $u=\varphi _{1}^{-1}(w),v=\varphi _{2}^{-1}(w).$ Then 
$\varphi _{1}(u)=\varphi _{1}(\varphi _{1}^{-1}(w))\in \lbrack u_{1},\infty
).$ Similarly $\varphi _{2}(v)\in \lbrack u_{1},\infty )$ and we finish as
above. In the case of nondegenerate Orlicz functions we simply get the
following equivalence:

$\left( a\right) $ there are $u_{0}>0$ and $C_{0}>0$ such that $C_{0}\varphi
_{1}^{-1}(u)\varphi _{2}^{-1}(u)\leq \varphi ^{-1}(u)$ for all $u\geq u_{0}.$

$\left( b\right) $ there are $u_{1}>0$ and $C_{1}>0$ such that $\varphi
(C_{1}uv)\leq \varphi _{1}(u)+\varphi _{2}(v)$ for all $u,v$ 

\hspace{5mm} with $u,v\geq u_{1}.$

\noindent
Finally, for all Orlicz functions, it is easy to show the following
equivalence:

$\left( a\right) $ for any $u_{0}>0$ there is $C_0 >0$ such that 
$C_{0}\varphi _{1}^{-1}(u)\varphi _{2}^{-1}(u)\leq \varphi ^{-1}(u)$ for all 
$u\leq u_{0}.$

$\left( b\right) $ for any $u_{1}>0$ there is $C_{1}>0$ such that $\varphi
(C_{1}uv)\leq \varphi _{1}(u)+\varphi _{2}(v)$ for all $u,v$ 

\hspace{5mm} with $\max \left\{ \varphi _{1}(u),\varphi _{2}(v)\right\} \leq u_{1}.$

\medskip
In the case $E = L^1$ the space $E_{\varphi}$ is an Orlicz space 
$L^{\varphi}$ and our Theorems 2--4 together with the equivalence of
conditions (\ref{left-ineq}) and (\ref{sufficiency}) give the following
results of Ando \cite{An60} and O'Neil \cite{ON65} (see also \cite{Ma89},
Theorems 10.2-10.4).

\vspace{3mm} 
Ando theorem (\cite{An60}, Theorem 1): \textit{Let $\mu$ be a non-atomic measure
and $0 < \mu(\Omega) < \infty$. For any $x \in L^{\varphi_1}$ and any $y \in
L^{\varphi_2}$ the product $x y \in L^{\varphi}$ if and only if there exist 
$C, u_0 > 0$ such that $\varphi(C uv) \leq \varphi_1(u) + \varphi_2(v)$ for
all $u, v \geq u_0$. }

\medskip 
O'Neil \cite[Theorem 6.5]{ON65} observed that the last condition on
Young functions is equivalent to the condition $\lim \sup_{u \rightarrow
\infty} \frac{\varphi _{1}^{-1} (u) \,\varphi _{2}^{-1}(u)}{\varphi ^{-1}(u)}
< \infty$. Krasnoselski{\u \i} and Ruticki{\u \i} \cite{KR61} noted that
relation on embedding of sets is equivalent to estimations of the norms,
that is, there is a number $A > 0$ such that $\| x y \|_{\varphi} \leq A
\|x\|_{\varphi_1} \| y \|_{\varphi_2}$ for all $x \in L^{\varphi_1}$ and $y
\in L^{\varphi_2}$.

\vspace{3mm} 
O'Neil theorems (\cite{ON65}, Theorems 6.6 and 6.7): \textit{(i) Let $\mu$ be
non-atomic measure with $\mu(\Omega) = \infty$. Then the following conditions are
equivalent: for any $x \in L^{\varphi_1}$ and any $y \in L^{\varphi_2}$ the
product $x y \in L^{\varphi}$ $\Longleftrightarrow $ there exist $C > 0$
such that $\varphi(C uv) \leq \varphi_1(u) + \varphi_2(v)$ for all $u, v > 0$
$\Longleftrightarrow $ $\sup_{u > 0} \frac{\varphi _{1}^{-1} (u) \,\varphi
_{2}^{-1}(u)}{\varphi ^{-1}(u)} < \infty \Longleftrightarrow $ there is a
number $B > 0$ such that $\| x y \|_{\varphi} \leq B \|x\|_{\varphi_1} \, \|
y \|_{\varphi_2}$ for all $x \in L^{\varphi_1}$ and $y \in L^{\varphi_2}$. }

\textit{(ii) Let $I=\mathbb{N}$ with the counting measure. The
following conditions are equivalent: for any $x\in l^{\varphi _{1}}$ and any 
$y\in l^{\varphi _{2}}$ the product $x\,y\in l^{\varphi }$ 
$\Longleftrightarrow $ there exist $C,u_{0}>0$ such that $\varphi (Cuv)\leq
\varphi _{1}(u)+\varphi _{2}(v)$ for all $0<u,v\leq u_{0}$ 
$\Longleftrightarrow $ $\lim \sup_{u\rightarrow 0^{+}}\frac{\varphi
_{1}^{-1}(u)\,\varphi _{2}^{-1}(u)}{\varphi ^{-1}(u)}<\infty
\Longleftrightarrow $ there is a number $D>0$ such
that $\Vert x\,y\Vert _{\varphi }\leq D\Vert x\Vert _{\varphi _{1}}\Vert
y\Vert _{\varphi _{2}}$ for all $x\in L^{\varphi _{1}}$ and $y\in L^{\varphi
_{2}}$.}

\vspace{3mm}
\begin{center}
\textbf{5. On the inclusion $M\left( E_{\varphi _{1}}, E_{\varphi }\right)
\hookrightarrow E_{\varphi _{2}} $}
\end{center}


We start by stating a crucial lemma, which in the case $E = L^1$ was proved in \cite{MN}.

\vspace{3mm} 
\textbf{Lemma 4.} \label{lem-3} \textit{If $\varphi \in \mathcal{Y}^{\left(
1\right) }\cup \mathcal{Y} ^{\left( 2\right) }$\ and $x = \sum_{k=1}^{N} c_{k} \chi_{A_k}, x \neq 0$ 
is a simple function, then $I_{\varphi }\left( \frac{x}{\| x \|_{E_{\varphi }}}\right)
=1.$}

\vspace{3mm} 
{\it Proof.} We follow arguments as it was done in the proof of
Lemma 3 in \cite{MN}. It is enough to show that the function 
\begin{equation*}
h( \lambda) = I_{\varphi }\left( \frac{x}{\lambda }\right) = \|
\sum_{k=1}^{N}\varphi \left( \frac{c_{k}}{\lambda }\right) \chi_{A_k} \|
_{E}
\end{equation*}
is continuous, non-increasing and $h: ( 0,c_{N}/a_{\varphi }) \rightarrow (
0, \infty)$. Suppose that $\varphi \in \mathcal{Y}^{\left( 1\right) }$. If 
$\lambda _{m} \rightarrow \lambda _{0}$, then 
\vspace{-2mm}
\begin{eqnarray*}
| h( \lambda _{m}) - h( \lambda _{0})| &\leq& \| \sum_{k=1}^{N} \left |
\varphi \left( \frac{c_{k}}{\lambda _{m}}\right) -\varphi \left( \frac{c_{k}
}{\lambda _{0}} \right) \right | \chi _{A_{k}} \|_{E} \\
&\leq& \sum_{k=1}^{N} \left | \varphi \left( \frac{c_{k}}{\lambda _{m}}
\right) -\varphi \left( \frac{c_{k}}{\lambda _{0}}\right) \right | \| \chi
_{A_{k}}\| _{E}\rightarrow 0
\end{eqnarray*}
as $m\rightarrow \infty$. Clearly, $h$ is non-increasing and 
\begin{equation*}
\lim_{\lambda \rightarrow 0^+}h( \lambda) \geq \lim_{\lambda \rightarrow
0^+} \varphi ( \frac{c_{1}}{\lambda }) \| \chi_{A_{1}} \| _{E} = \infty, ~
\lim_{\lambda \rightarrow c_{N}/a_{\varphi }}h\left( \lambda \right) =0,
\end{equation*}
since for $\lambda >\frac{c_{N-1}}{a_{\varphi }}$ we have that $h\left( \lambda
\right) =\varphi \left( \frac{c_{N}}{\lambda }\right) \left\Vert \chi
_{A_{N}}\right\Vert _{E}.$ Consequently, there is a number $\lambda _{0}\in
\left( 0,c_{N}/a_{\varphi }\right) $ with $I_{\varphi }\left( \frac{x}{
\lambda _{0}}\right) =1.$

If $\varphi \in \mathcal{Y}^{\left( 2\right) }$ the proof is the same as in 
\cite{MN} and Lemma 4 is proved. \qed

\vspace{3mm} 
\noindent 
The following result is a generalization of Theorem 1 from \cite{MN}.

\vspace{3mm} 
\textbf{THEOREM 5}. \label{thm-4} {\it Suppose $E$ is a Banach ideal
space with the Fatou property and ${\it supp} E=\Omega$. Let
$\varphi, \varphi _{1}, \varphi _{2}$ be Young functions. Assume that 
the condition $\varphi ^{-1}\prec \varphi_{1}^{-1}\varphi _{2}^{-1}$ holds: 

\vspace{-3mm}

\begin{itemize}
\item[$(i)$] for all arguments.

\item[$(ii)$] for large arguments and $L^{\infty} \hookrightarrow E$.

\item[$(iii)$] for small arguments and $E\hookrightarrow L^{\infty }$. 
\end{itemize}
\vspace{-2mm}

\noindent
Then $M\left( E_{\varphi _{1}}, E_{\varphi }\right)
\hookrightarrow E_{\varphi_{2}}.$}

\vspace{3mm}

{\it Proof.} We apply the technique from the proof of Theorem 1 in \cite{MN}.
The case $(i)$ follows in the same way as in \cite{MN} with one restriction.
Namely, if $x\in M( E_{\varphi _{1}}, E_{\varphi })$ is a simple function,
then $x$ need not belong to $E_{\varphi _{2}}.$ Consider the case 
$x=\sum_{i=1}^{n}a_{i}\chi _{A_{i}}$ and $\chi _{A_{i}}\notin E$ for some $i.$
Then there is an increasing sequence $\left( A_{i}^{k}\right)
_{k=1}^{\infty} $ of subsets of $A_{i}$ satisfying $\bigcup_{k=1}^{\infty}
A_{i}^{k} = A_{i}$ and $\chi _{A_{i}^{k}} \in E$ for each $k.$ Taking 
$x_{k}=\sum_{i=1}^{n}a_{i}\chi _{A_{i}^{k}}$ we get $x_{k}\in E_{\varphi
_{2}} $ for each $k.$ Then we follow just the proof of Theorem 1 in \cite{MN}
and get $\| x_{k} \| _{M(E_{\varphi _{1}}, E_{\varphi }) }\geq \frac{1}{D}
\| x_{k} \| _{\varphi _{2}}$ for each $k.$ But $0\leq x_{k}\leq x$ and 
$x_{k}\uparrow x$. Since $E_{\varphi _{2}}$ has the Fatou property it follows
that $\| x \|_{M\left( E_{\varphi _{1}}, E_{\varphi }\right) }\geq 
\frac{1}{D} \| x \|_{\varphi _{2}}$ for each simple function $x \in M(
E_{\varphi _{1}}, E_{\varphi })$.

\medskip
$(ii)$ Assume that $L^{\infty }\hookrightarrow E.$ Take $\alpha >a_{\varphi
_{2}}$ with $\varphi _{2}(\alpha )\Vert \chi _{\Omega }\Vert _{E}<\frac{1}{2}$.
Applying Lemma 3 we find a contant $D_{1}\geq D$ such that 
\vspace{-2mm}
\begin{equation*}
\varphi ^{-1}(u)\leq D_{1}\varphi _{1}^{-1}(u)\varphi _{2}^{-1}(u)~~\mathrm{
for~any}~u\geq \varphi _{2}(\alpha ).
\end{equation*}
Observe that 
\begin{equation}
I_{\varphi _{2}}(z\chi _{B})\geq 1/2  \label{ineq-15}
\end{equation}
for any $z = z(t)$ with $I_{\varphi _{2}}\left( z\right) =1$, where $B=\{t\in 
\mathrm{supp}\,z: |z(t)|\geq \alpha \}$. Really, otherwise 
\vspace{-2mm}
\begin{equation*}
1=I_{\varphi _{2}}(z)\leq I_{\varphi _{2}}(z\chi _{B})+I_{\varphi
_{2}}(z\chi _{\Omega \backslash B})<\frac{1}{2}+\varphi _{2}(\alpha )\Vert
\chi _{\Omega }\Vert _{E}<1,
\end{equation*}
and we get a contradiction. Although some steps are similar as in the proof
of Theorem 1 from \cite{MN} we present the whole proof for the sake of
completeness. Assume that $\varphi $ and $\varphi _{2}$ are in $\mathcal{Y}
^{(1)}\cup \mathcal{Y}^{(2)}.$ Let $x\in M(E_{\varphi _{1}}, E_{\varphi })$.
First suppose that $x$ is a simple function. Then $x\in E_{\varphi _{2}},$
because $\varphi _{2}\circ (\lambda |x|)$ is a simple function in $E$ for
some $\lambda $ and $\chi _{A}\in L^{\infty }\hookrightarrow E$ for each 
$A\in \Sigma $. Consequently, 
\begin{equation*}
y(t)=\varphi _{2}\left( \frac{|x(t)|}{\Vert x\Vert _{E_{\varphi _{2}}}}
\right) <\infty ~~\mathrm{for}~~\mu -\mathrm{a.e.}~t\in \Omega .
\end{equation*}
Set 
\begin{equation*}
z(t)=\left\{ 
\begin{array}{ccc}
\varphi _{1}^{-1}(y(t)) & \text{if} & 0<y(t)<\infty , \\ 
0 & \text{if} & y(t)=0.
\end{array}
\right.
\end{equation*}
Then $I_{\varphi _1}(z)\leq I_{\varphi _2}(\frac{x}{\| x\|
_{E_{\varphi _2}}})\leq 1$ implies $\| z\|_{E_{\varphi _1}}\leq 1$
and, by the assumption, we have $zx\in E_{\varphi }$. Denote 
\begin{equation*}
A=\{t\in \mathrm{supp}\,y:\frac{|x(t)|}{\| x\|_{E_{\varphi _{2}}}}
<\alpha \}\text{ and } B=\{t\in \mathrm{supp}\,y:\frac{|x(t)|}{\| x\|
_{E_{\varphi _{2}}}}\geq \alpha \}.
\end{equation*}
Then, for $\mu $-a.e. $t\in B$, 
\begin{equation*}
z(t)\,\frac{|x(t)|}{\Vert x\Vert _{E_{\varphi _{2}}}}=\varphi
_{1}^{-1}(y(t))\varphi _{2}^{-1}(y(t))\geq \frac{1}{D_{1}}\varphi
^{-1}(y(t)),
\end{equation*}
so that
\begin{equation*}
\varphi \left( D_{1}z(t)\frac{2|x(t)|}{\| x\|_{E_{\varphi _{2}}}}
\right) \geq \varphi (2\varphi ^{-1}(y(t)))\geq 2\varphi (\varphi
^{-1}(y(t)))=2y(t),
\end{equation*}
where the last equality follows from the fact that $\varphi \in \mathcal{Y}
^{(1)}\cup \mathcal{Y}^{(2)}$. By Lemma 4, we have $I_{\varphi _{2}}\left( 
\frac{x}{\| x\| _{E_{\varphi _{2}}}}\right) =1$ and, consequently, by (\ref{ineq-15}), 
\vspace{-2mm}
\begin{equation*}
I_{\varphi }\left( D_{1}z\,\frac{2x}{\Vert x\Vert _{E_{\varphi _{2}}}}\chi
_{B}\right) \geq 2I_{\varphi _{2}}\left( \frac{x}{\Vert x\Vert _{E_{\varphi
_{2}}}}\chi _{B}\right) \geq 1.
\end{equation*}
Thus $\| zx\|_{E_{\varphi }}\geq \frac{1}{2D_{1}}\| x\|_{E_{\varphi _{2}}}$ and so 
$\| x\|_{M\left( E_{\varphi_{1}}, E_{\varphi }\right) }\geq \frac{1}{2D_{1}}\| x\|_{E_{\varphi
_{2}}}$. For arbitrary function $x\in M\left( E_{\varphi _{1}}, E_{\varphi
}\right) $ we can follow the proof of Theorem 1 from \cite{MN},
because the space $E_{\varphi _{2}}$ has the Fatou property provided that $E$ has it.

If $\varphi $ or $\varphi _{2}$ is in $\mathcal{Y}^{(3) }$ we follow again
the same way as in case 2 of the proof of Theorem 1 in \cite{MN} to show
that $\| x \|_{M\left( E_{\varphi _{1}}, E_{\varphi }\right) }\geq \frac{1}{2D_{1}}
 \, \| x \| _{E_{\varphi _{2}}}$.

$(iii)$ Let $E\overset{A}{\hookrightarrow }L^{\infty }.$ First observe that
if $\| u \|_{E_{\varphi }} \leq 1$, then ${\rm ess\, sup}_{t \in \Omega} \,
| \varphi (u(t))| \leq A$. Furthermore, by assumption $(iii)$ and Lemma 3, 
there exists a contant $D_{2}$ such that 
\begin{equation*}
\varphi ^{-1}(u)\, \leq D_{2} \varphi _{1}^{-1} (u) \varphi _{2}^{-1} (u) ~ 
\mathrm{for ~ any} ~ 0 < u\leq A.
\end{equation*}
The rest of the proof goes as in the case $(i)$ (see also the proof of
Theorem 1 in \cite{MN}). The proof is complete. \qed

\medskip

To find cases when the condition $\varphi ^{-1}\prec \varphi_{1}^{-1}\varphi_{2}^{-1}$ 
is necessary for the imbedding 
$M\left(E_{\varphi_{1}}, E_{\varphi }\right) \hookrightarrow E_{\varphi _{2}}$ we
will take as $E$ a symmetric function space on $I$. Then $E_{\varphi}$ is
also a symmetric function space and an easy calculation gives that the
fundamental function $f_{E_{\varphi}}$ of $E_{\varphi}$ is equal to 
$f_{E_{\varphi}}(t) = \frac{1}{\varphi^{-1}(1/f_E(t))} ~~ \mathrm{for} ~~ t
\in (0, m(I)).$

\vspace{3mm}
{\bf THEOREM 6}. \label{thm5} {\it Let $E$ be a symmetric function
space on $I$ and let $\varphi, \varphi _{1}, \varphi _{2}$ be Young
functions. Suppose 
\begin{equation}
M\left( E_{\varphi _{1}}, E_{\varphi }\right) \hookrightarrow E_{\varphi
_{2}}.  \label{fund 1}
\end{equation}

\vspace{-2mm}

\begin{itemize}
\item[$(i)$] If there are numbers $a,b>0$ such that $\frac{
f_{E_{\varphi }} (t) }{f_{E_{\varphi _{1}}}( t) \, t^{a}} = \frac{
\varphi_{1}^{-1} (1/f_{E}(t)) }{\varphi^{-1} (1/f_{E}(t)) \, t^{a}}$ is a
non-decreasing function of $t$ on an interval $(0, b)$ and $E_{a}\not= \{ 0
\}$, then $\varphi ^{-1} \prec \varphi _{1}^{-1}\varphi _{2}^{-1}$ for large
arguments. 

\vspace{-3mm} 

\item[$(ii)$] Let $b_{\varphi} = \infty$. If there is a number $a>0$ such that 
$\frac{f_{E_{\varphi }} (t)}{f_{E_{\varphi _{1}}} ( t) \, t^{a}} = \frac{
\varphi_{1}^{-1} (1/f_{E}(t)) }{\varphi^{-1} (1/f_{E}(t)) \, t^{a}}$ is a
non-decreasing function of $t$ on $(0, \infty), L^{\infty
}\not\hookrightarrow E$ and ${\it supp} E_{a} = {\it supp} E$, then $\varphi^{-1}
\prec \varphi_1^{-1} \varphi_2^{-1}$ for all arguments. 
\end{itemize}
}

{\it Proof.} (i) Assume $f_{E_{\varphi }}( 0^+) =0$ and suppose that the condition 
$\varphi ^{-1}\prec \varphi_{1}^{-1}\varphi _{2}^{-1}$ is not satisfied for large arguments, 
i.e., there is a sequence $\left( u_{n}\right) $ tending to infinity such that for any 
$n\in \mathbb{N}$ 
\begin{equation*}
2^{n}\varphi _{1}^{-1}(u_{n})\varphi _{2}^{-1}(u_{n})\leq \varphi
^{-1}(u_{n}).
\end{equation*}
It is enough to find a sequence $\left( x_{n}\right) $ both in $M(E_{\varphi
_{1}}, E_{\varphi })$ and $E_{\varphi _{2}}$ such that 
\begin{equation*}
\frac{\| x_{n}\|_{E_{\varphi _2}}}{\| x_{n}\|_{M(E_{\varphi_{1}}, E_{\varphi })}}\longrightarrow \infty .
\end{equation*}
Analogously as in Theorem 3(i) for each $u_{n}$ one can find measurable set 
$A_{n}$ satisfying $\| u_{n}\chi _{A_{n}}\| _{E}=1$. Define 
\begin{equation*}
x_{n}=\varphi_{2}^{-1}(u_{n})\chi _{A_{n}}.
\end{equation*}
Then $\| x_{n} \|_{E_{\varphi _{2}}} = 1$. In fact, if $\varphi _{2}\in \mathcal{Y}^{(1)}\cup \mathcal{Y}^{(2)}$,
then 
\begin{equation*}
I_{\varphi _{2}}( x_{n}) = \| \varphi _{2}\left( \varphi
_{2}^{-1}\left( u_{n}\right) \right) \chi _{A_{n}} \|_{E}=u_{n}\left\Vert \chi _{A_{n}}\right\Vert _{E}=1.
\end{equation*}
If $\varphi _{2}\in \mathcal{Y}^{(3)}$, then there is $N_{0}$ with $\varphi
_{2}^{-1}\left( u_{n}\right) =b_{\varphi _{2}}$ for $n\geq N_{0},$ whence 
$I_{\varphi _{2}}\left( x_{n}\right) \leq 1$ and $I_{\varphi _{2}}\left(
x_{n}/\lambda \right) =\infty $ for $0<\lambda <1.$ Thus,  
$\| x_{n}\|_{E_{\varphi _{2}}} = 1$ for sufficiently large $n.$
\vspace{1mm}

Putting $t_{n}=m(A_{n})$ we obtain by symmetry of $E$ that 
$f_{E}(t_{n})=\Vert \chi _{\lbrack 0,t_{n}]}\Vert _{E}=\Vert \chi
_{A_{n}}\Vert _{E}=\frac{1}{u_{n}}\rightarrow 0$ as $n\rightarrow \infty $
and so $t_{n}\rightarrow 0$ as $n\rightarrow \infty $. Therefore, according
to Theorem 1(iii), for $t_{n}\in (0,b)$, we obtain 
\vspace{-2mm}
\begin{eqnarray*}
\| x_{n} \|_{M(E_{\varphi _{1}}, E_{\varphi })} &=&\varphi
_{2}^{-1}(u_{n})\,\Vert \chi _{A_{n}}\Vert _{M(E_{\varphi _{1}}, E_{\varphi
})}=\varphi _{2}^{-1}(u_{n})f_{M(E_{\varphi _{1}}, E_{\varphi })}(t_{n}) \\
&\leq &
2 \, \frac{1}{a}\,\varphi _{2}^{-1}(u_{n})\frac{f_{E_{\varphi }}(t_{n})}{
f_{E_{\varphi _{1}}}(t_{n})} \leq 2\, \frac{\varphi ^{-1}(u_{n})}{
a\,2^{n}\,\varphi _{1}^{-1}(u_{n})}\frac{f_{E_{\varphi }}(t_{n})}{
f_{E_{\varphi _{1}}}(t_{n})} \\
&=&
\frac{2}{a\,2^{n}}\,\frac{\varphi ^{-1}(u_{n})}{\varphi _{1}^{-1}(u_{n})}
\,\frac{\varphi _{1}^{-1}(1/f_{E}(t_{n}))}{\varphi ^{-1}(1/f_{E}(t_{n}))}=
\frac{2}{a\,2^{n}}\rightarrow 0~\mathrm{as}~n\rightarrow \infty ,
\end{eqnarray*}
which finishes the proof.

In the case when $f_{E_{\varphi }}( 0^+) > 0$ we have $f_E( 0^+) > 0$ which implies 
$b_{\varphi} < \infty$ and estimate on Young functions is automatically 
satisfied. 

$(ii)$ This part is analogous to the above and Theorem 3(ii). \qed

\medskip 
\noindent
Combining Theorems 5 and 6, we obtain the following result:

\vspace{3mm} 
{\bf Corollary 3}. \label{cor-3} {\it Let $E$ be a symmetric function space 
on I with the Fatou property. Let $\varphi, \varphi _{1},
\varphi _{2}$ be Young functions. 
\vspace{-3mm} 

\begin{itemize}
\item[$(i)$] Suppose $L^{\infty }\hookrightarrow E, E_{a}\not= \{ 0
\} $ and that there are numbers $a, b>0$ such that $\frac{ \varphi _{1}^{-1} (
1/f_{E}(t)) }{\varphi ^{-1}( 1/f_{E}(t)) \, t^{a}}$ is a non-decreasing
function of $t$ on the interval $(0, b)$. Then 
$M(E_{\varphi_{1}}, E_{\varphi }) \hookrightarrow E_{\varphi _{2}}$ if and only if $\varphi
^{-1}\prec \varphi _{1}^{-1}\varphi _{2}^{-1}$ for large arguments. 
\vspace{-2mm} 

\item[$(ii)$] Assume $b_{\varphi} = \infty$, $L^{\infty }\not\hookrightarrow E$, ${\it supp} E_{a} = I$ 
and that there is a number $a>0$ such that $\frac{\varphi_{1}^{-1}
( 1/f_{E}(t) ) }{\varphi ^{-1}(1/f_{E}(t)) \,t^{a}}$ is a non-decreasing
function on $(0, \infty)$. Then $M(E_{\varphi _{1}}, E_{\varphi })
\hookrightarrow E_{\varphi _{2}}$ if and only if $\varphi ^{-1}\prec
\varphi_{1}^{-1}\varphi _{2}^{-1}$ for all arguments. 
\end{itemize}
}

Now, if we take in Corollary 3 as $E=L^{1}$ we obtain the results, which 
give an answer for the problem posed in the book \cite{Ma89} (Problem 4, p. 77) in
the case of Orlicz spaces (under an additional assumption):

(i) Let $I = [0, 1]$ and let $\frac{ \varphi _{1}^{-1} (u)\, u^a }{\varphi
^{-1}(u)}$ be a non-increasing function for some $a > 0$ and suffieciently large u. 
Then $M\left(L^{\varphi_1}, L^{\varphi }\right) \hookrightarrow L^{\varphi _2}$ if and
only if $\varphi ^{-1}\prec \varphi _{1}^{-1}\varphi _{2}^{-1}$ for large
arguments. 

(ii) Let $I = [0, \infty)$ and let $\frac{ \varphi _{1}^{-1} (u)\, u^a }{\varphi ^{-1}(u)}$ 
be a non-increasing function for some $a > 0$ and all $u > 0$. Then 
$M\left( L^{\varphi_1}, L^{\varphi }\right) \hookrightarrow L^{\varphi _2}$
if and only if $\varphi ^{-1}\prec \varphi _{1}^{-1}\varphi _{2}^{-1}$ for
all arguments.

The monotonicity assumption in (i) is essential for the equivalence (see Example 9 (g) below).

\vspace{3mm}
\begin{center}
\textbf{6. On the equality $M\left( E_{\varphi _{1}}, E_{\varphi }\right) =
E_{\varphi _{2}} $}
\end{center}


Putting together Theorems 2 and 5 we obtain sufficient conditions for coincidence of the
space of pointwise multipliers $M\left( E_{\varphi _{1}}, E_{\varphi }\right)$
with $E_{\varphi _{2}}$.

\vspace{3mm} 
{\bf Corollary 4}. \label{cor-4} {\it Let $\varphi, \varphi _{1}$ and $\varphi _{2}$ be Young 
functions. Suppose $E$ is a Banach ideal space with the Fatou property and 
${\it supp} E =\Omega$. Assume also that at least one of the following conditions holds: 
\vspace{-3mm} 

\begin{itemize}
\item[$(i)$] $\varphi _{1}^{-1}\varphi _{2}^{-1}\approx \varphi
^{-1} $ for all arguments. 
\vspace{-2mm} 

\item[$(ii)$] $\varphi _{1}^{-1}\varphi _{2}^{-1}\approx \varphi
^{-1}$ for large arguments and $L^{\infty }\hookrightarrow E$ . 
\vspace{-2mm} 

\item[$(iii)$] $\varphi _{1}^{-1}\varphi _{2}^{-1}\approx \varphi
^{-1} $ for small arguments and $E\hookrightarrow L^{\infty }.$ 
\end{itemize}
\vspace{-2mm}
Then $M\left( E_{\varphi _{1}}, E_{\varphi }\right) = E_{\varphi _{2}}$ 
with equivalent norms.
}

\vspace{3mm}
Taking into account Corollary 1 and 3, we obtain

\vspace{3mm}
{\bf Corollary 5}. \label{cor-5} {\it Let $E$ be a symmetric function space 
with the Fatou property. Let $\varphi, \varphi _{1}, \varphi
_{2}$ be Young functions. 
\vspace{-3mm} 

\begin{itemize}
\item[$(i)$] Suppose $L^{\infty }\hookrightarrow E, E_{a}\not= \{
0\} $ and that there are numbers $a, b>0$ such that $\frac{\varphi _{1}^{-1} (
1/f_{E}(t) }{\varphi ^{-1} (1/f_{E}(t))\, t^a}$ is a non-decreasing function
of $t$ on the interval $(0, b)$. Then $M( E_{\varphi_1}, E_{\varphi }) \newline =
E_{\varphi _2}$ if and only if $\varphi ^{-1}\approx \varphi
_{1}^{-1}\varphi _{2}^{-1}$ for large arguments.
\vspace{-2mm} 

\item[$(ii)$] Assume $b_{\varphi} = \infty$, $L^{\infty }\not\hookrightarrow E$, 
${\it supp} E_{a} = {\it supp} E$ and that there is  a number $a>0$ such that $\frac{\varphi
_{1}^{-1}(1/f_{E}(t))}{\varphi ^{-1}(1/f_{E}(t))\,t^{a}}$ is a non-decreasing
function on $(0,\infty )$. Then $M(E_{\varphi _{1}}, E_{\varphi })=E_{\varphi
_{2}}$ if and only if $\varphi ^{-1}\approx \varphi _{1}^{-1}\varphi
_{2}^{-1}$ for all arguments. 
\end{itemize}
}

\begin{center}
\textbf{7. On the construction of a Young function generating the space 
$M\left(E_{\varphi _{1}}, E_{\varphi }\right) $ }
\end{center}

The following questions arises: having two Young functions 
$\varphi_{1},\varphi $ how can one find a Young function $\varphi _{2}$
satisfying $\varphi _{1}^{-1}\varphi _{2}^{-1}\approx \varphi ^{-1}$? Does
such a function always exist?

It appears that such a function may not exist. The following example
describes such possibility.

\vspace{3mm}
\textbf{Example 5}. Let $\varphi (u)=u^{2},\varphi _{1}(u)=u$ and 
$E=L_{t}^{\infty }[0,1]$ with the norm $\| x\| _{E} =$ ${\rm ess\,sup}
_{t\in \left[ 0,1\right] } | t\,x(t) |$. The equivalence 
$\varphi _{1}^{-1}\varphi _{2}^{-1}\approx \varphi ^{-1}$ means that $u\,\varphi
_{2}^{-1}\approx \sqrt{u}$, i.e., $\varphi _{2}^{-1}\approx 1/\sqrt{u}$, 
which is not possible for any Young function $\varphi _{2}$. Moreover, 
$M(E_{\varphi _{1}}, E_{\varphi })$ is not a Calder\'{o}n-Lozanovski\u{\i}
space of the form $E_{\varphi _{3}}$ for any Young function $\varphi _{3}$.
In fact, $E_{\varphi }=L_{\sqrt{t}}^{\infty }[0,1]$ and 
\begin{equation*}
M(E_{\varphi _{1}}, E_{\varphi }) = M(L_{t}^{\infty }[0,1], L_{\sqrt{t}}^{\infty
}[0,1]) = L_{1/\sqrt{t}}^{\infty }[0,1].
\end{equation*}
This space cannot be of the form $E_{\varphi _{3}}$ since $\chi _{\lbrack
0,1]}\in E_{\varphi _{3}}$ and $\Vert \chi _{\lbrack 0,1]}\Vert _{E_{\varphi
_{3}}}=1/\varphi _{3}^{-1}(1)$, but $\chi _{\lbrack 0,1]}\not\in L_{1/\sqrt{t
}}^{\infty }[0,1].$

\medskip 
We have seen in the proof of Theorem 2 that the inequality $\varphi
_{1}^{-1}(u)\,\varphi _{2}^{-1}(u)\leq \varphi ^{-1}(u)$ for all $u>0$ gives 
that $\varphi (uv)\leq \varphi _{1}(v)+\varphi _{2}(u)$ for all $u,v>0$ and 
that the last estimate suggests to consider an operation on two Young functions 
$\varphi _{1},\varphi $ and compare it with $\varphi _{2}$. Define a new
function $\varphi \ominus \varphi _{1}:[0,\infty )\rightarrow \lbrack
0,\infty ]$ by the formula 
\vspace{-1mm}
\begin{equation*}
\left( \varphi \ominus \varphi _{1}\right) \left( u\right) =\sup_{v\geq
0}\left[ \varphi \left( uv\right) -\varphi _{1}\left( v\right) \right] .
\end{equation*}

\vspace{-2mm}

\noindent
We may say that $\varphi \ominus \varphi _{1}$ is the conjugate
(complementary) function (in the sense of Young) to $\varphi _{1}$ with respect
to $\varphi $. In particular, if $\varphi (u)=u$, then $\varphi \ominus
\varphi _{1}=\varphi _{1}^{\ast }$\ is the usual conjugate (complementary)
function (in sense of Young) to $\varphi _{1}.$ This operation on the class of
N-functions was defined by Ando \cite[p. 180]{An60} and on the class of
extended Young functions (by word "extended Young" functions we mean
nondecreasing convex functions $\varphi :[0,\infty )\rightarrow \lbrack
0,\infty ]$ with $\varphi (0)=0$ and they can be trivial) by O'Neil \cite[p.
325]{ON65} and he referred to Ando.

Note that it can happen that the function $\varphi \ominus \varphi
_{1}\left( u\right) =\infty $ for $u>0$, and then we have that the
corresponding Orlicz space is the zero space. To avoid a confusion when 
$\max \{b_{\varphi }, b_{\varphi _{1}} \} <\infty $ since then we will have symbol 
$\infty - \infty$ we better skip this case.

Moreover, in the case when we work with sequence spaces (or in case 
$E\hookrightarrow L^{\infty }$) it is reasonable to define $\varphi \ominus \varphi _{1}$ 
in a different way, namely,  
\begin{equation*}
\left( \varphi \ominus \varphi _{1}\right) _{0}\left( u\right) =\sup_{0\leq
v\leq 1}\left[ \varphi \left( uv\right) -\varphi _{1}\left( v\right)
\right],
\end{equation*}

\vspace{-2mm}

\noindent
since then only the behaviour of the functions in a neighbourhood of zero
is important. Djakov and Ramanujan \cite{DR00} proved that in the case of Orlicz 
sequence spaces we have that $M(l^{\varphi_1}, l^{\varphi}) =
l^{\varphi_2}$, where $\varphi_2 = (\varphi \ominus \varphi _{1})_0$. It is easy to see that
the function $(\varphi \ominus \varphi _1) _{0}$ is smaller than $\varphi \ominus \varphi _1$ 
and it can be different from $\varphi \ominus \varphi _1 $. 

\vspace{3mm} 
\textbf{Example 6}. Let $\varphi (u) = u^{p}/p, \varphi _1(u) = u^{p_1}/p_1$
with $1 \leq p, p_1 < \infty$. If $p > p_1$, then $(\varphi \ominus \varphi
_1) (u) = \infty$ for $u > 0$ and 
\begin{equation*}
~~(\varphi \ominus \varphi _1) _{0}(u) = 
\begin{cases}
0 & ~\mathrm{if} ~0 \leq u \leq (p/p_1)^{1/p}, \\ 
\frac{u^p}{p} - \frac{1}{p_1} & ~\mathrm{if} ~u \geq (p/p_1)^{1/p}.
\end{cases}
\end{equation*}
If $p = p_1$, then 
\begin{equation*}
~~(\varphi \ominus \varphi _1)(u) = 
\begin{cases}
0 & ~\mathrm{if} ~0 \leq u \leq 1, \\ 
\infty & ~\mathrm{if} ~u > 1,
\end{cases}
\end{equation*}
and 
\begin{equation*}
~~(\varphi \ominus \varphi _1) _{0}(u) = 
\begin{cases}
0 & ~\mathrm{if} ~0 \leq u \leq 1, \\ 
\frac{u^p - 1}{p} & ~\mathrm{if} ~u \geq 1.
\end{cases}
\end{equation*}
If $p < p_1$, then $(\varphi \ominus \varphi _1)(u) = \frac{u^{p_2}}{p_2}$,
where $\frac{1}{p_2} = \frac{1}{p} - \frac{1}{p_1}$ and 
\begin{equation*}
~~(\varphi \ominus \varphi _1) _{0}(u) = 
\begin{cases}
\frac{u^{p_2}}{p_2} & ~\mathrm{if} ~0 \leq u \leq 1, \\ 
\frac{u^p}{p} - \frac{1}{p_1} & ~\mathrm{if} ~u \geq 1.
\end{cases}
\end{equation*}

\vspace{3mm} 
\textbf{Example 7}. Let 
\begin{equation*}
~~\varphi (u) = 
\begin{cases}
0 & ~\mathrm{if} ~0 \leq u \leq 1, \\ 
u-1 & ~\mathrm{if} ~u \geq 1,
\end{cases}
\end{equation*}
and $\varphi_1(u) = u^2$. Then 
\begin{equation*}
~~ \varphi_2(u) = (\varphi \ominus \varphi _1)(u) = 
\begin{cases}
0 & ~\mathrm{if} ~0 \leq u \leq 2, \\ 
\frac{u^2}{4} - 1 & ~\mathrm{if} ~u \geq 2,
\end{cases}
\end{equation*}
\begin{equation*}
~~ \varphi_3(u) = (\varphi \ominus \varphi _2)(u) = 
\begin{cases}
0 & ~\mathrm{if} ~0 \leq u \leq 1/2, \\ 
2u - 1 & ~\mathrm{if} ~1/2 \leq u \leq 1, \\ 
u^2 & ~\mathrm{if} ~u \geq 1,
\end{cases}
\end{equation*}
and $\varphi_4(u) = (\varphi \ominus \varphi _3)(u) = \varphi _2(u)$ for all 
$u \geq 0$. The last equality was proved by O'Neil \cite[p. 325]{ON65}. 
For Orlicz spaces considered on $I = [0, \infty)$ we have 
\begin{equation} \label{equa24}
M(L^{\varphi_1}, L^{\varphi}) = M(L^2, L^1+L^{\infty}) = L^2 + L^{\infty} = L^{\varphi_2}
\end{equation}
and
\begin{equation} \label{equa25}
M(L^{\varphi_2}, L^{\varphi}) = M(L^2 + L^{\infty}, L^1+L^{\infty}) = L^{\varphi_3} = L^2 + L^{\infty} . 
\end{equation}
The second equality in (\ref{equa24}) we can get in the following way: if $x \in L^2+L^{\infty}$ and $y \in L^2$, 
then
\begin{eqnarray*}
\| x y\|_{L^1+L^{\infty}} 
&=& 
\int_0^1 (xy)^{\ast}(t) dt \leq  (\int_0^1 x^{\ast}(t/2)^2 dt)^{1/2} \, (\int_0^1 y^{\ast}(t/2)^2 dt)^{1/2} \\
&\leq&
2 \|x\|_{L^2+L^{\infty}} \| y\|_{L^2},
\end{eqnarray*}
and, hence, $L^2+L^{\infty} \overset{2}{\hookrightarrow } M(L^2, L^1+L^{\infty})$. This embedding also 
follows from Theorem 2(i) since $\varphi_1^{-1}(u) \varphi_2^{-1}(u) \leq 2 \varphi^{-1}(u)$ for all $u > 0$.

On the other hand, the function $y_t(s) = \chi_{[0, t]}(s)/\max(1, t) \in L^1 \cap L^{\infty}$ and 
$\| y_t\|_{ L^1 \cap L^{\infty}} = 1$. Thus, by the general property in (vii) and Theorem 1(i),
\begin{eqnarray*}
\| x\|_{M(L^2, L^1+L^{\infty})} 
&=& 
\| x \|_{M(L^1\cap L^{\infty}, L^2)} = \| x^{\ast} \|_{M(L^1\cap L^{\infty}, L^2)} \\
&\geq&
\| x^{\ast} y_t\|_{L^2} = \frac{1}{\max(1, t)}\| x^{\ast} \chi_{[0, t]}\|_{L^2} ~~ {\rm for ~any} ~t > 0.
\end{eqnarray*}
Hence,
\vspace{-2mm}
\begin{eqnarray*}
\| x \|_{M(L^2, L^1+L^{\infty})}
& \geq& 
\sup_{t > 0} \frac{1}{\max(1, t)}\| x^{\ast} \chi_{[0, t]}\|_{L^2} = (\int_0^1 x^{\ast}(s)^2 ds )^{1/2} \\
&\geq& 
\frac{1}{\sqrt{2}} \|x \|_{L^2+L^{\infty}},
\end{eqnarray*}

\vspace{-4mm}

\noindent
and we have the reverse embedding $M(L^2, L^1+L^{\infty}) \overset{\sqrt{2}}{\hookrightarrow } L^2+L^{\infty}$. 
This embedding does not follow from Theorem 5(i) or Corollary 4(i) since 
$\lim_{u \rightarrow 0^+} \frac{ \varphi^{-1}(u)}{\varphi_1^{-1}(u) \varphi_2^{-1}(u)} = 
\lim_{u \rightarrow 0^+} \frac{u+1}{\sqrt{u} \cdot 2 \sqrt{u+1}} = \infty$. This is also 
not a contradiction with Theorem 6(ii), Corollary 3(ii) and Corollary 5(ii) since the function 
$\varphi_1^{-1}(1/t)/[ \varphi^{-1}(1/t) t^a]$ is not non-decreasing for any $a > 0$.

The second equality in (\ref{equa25}) follows from Corollary 4(i) since 
$\varphi^{-1}(u) \leq \varphi_2^{-1}(u) \varphi_3^{-1}(u) \leq 2 \varphi^{-1}(u)$ for all $u > 0$.

\vspace{3mm} 
Some properties of operation $\varphi \ominus \varphi _1$ are collected
in the next lemma (part (iii) in the Lemma 5 below was proved in \cite[Theorem 3]{ZR67} with
some additional assumptions; cf. also \cite{Ma89} and \cite{MP89}).

\vspace{3mm}
{\bf Lemma 5.} \label{lem-4} {\it Let $\ \varphi, \varphi_1 $ be two
Young functions with $\max \{b_{\varphi }, b_{\varphi _{1}} \} = \infty $
and $\varphi_2 = \varphi \ominus \varphi _1$. 

\vspace{-3mm} 

\begin{itemize}
\item[$(i)$] The function $\varphi _{2}$ is non-decreasing, convex, left-continuous 
on $[0,\infty )$ with $\varphi _{2}(0)=0$ and it can be $\infty $ on $(0,\infty )$. 
\vspace{-2mm} 

\item[$(ii)$] We have 
\begin{equation*}
\varphi _{1}^{-1} (u) \, \varphi _{2}^{-1} (u) \leq 2 \, \varphi ^{-1} (u) ~ 
{\ for ~ all } ~ u\geq 0.
\end{equation*}

\vspace{-2mm} 

\item[$(iii)$] If $b_{\varphi } = b_{\varphi _1}=\infty $ (which
means that $\varphi, \varphi_1$ are, in fact, Orlicz functions) and for any $v
> 0$ the function $\frac{\varphi_1(u)}{\varphi(u v)}$ is equivalent to a
non-decreasing function, then $\varphi _{1}^{-1}\varphi _{2}^{-1}\approx
\varphi ^{-1}$ for all arguments. 
\end{itemize}
}

{\bf Remark 4}. This lemma gives a constructive way to define the function $\varphi_2$ such that 
$M(E_{\varphi_1}, E_{\varphi}) \hookrightarrow E_{\varphi_2}$ (see Theorem 5) and 
$E_{\varphi_2} \hookrightarrow M(E_{\varphi_1}, E_{\varphi})$ (see Theorem 2).
\medskip

\textit{Proof of Lemma 5.} (i) Of course, $\varphi _2(0) =0$ and $\varphi _2$ is
a non-decreasing function together with $\varphi $. Moreover, $\varphi _2$ is a convex
function since $\varphi$ is convex. We only need to show that $\varphi_2$ is left-continuous 
at $u_0 > 0$. We consider two cases.
\smallskip

$1^{\small 0}.$ Let $0<\varphi_2(u_0) < \infty $. Suppose, on the contrary, that $\varphi _2$ is 
not left-continuous at $u_0$. Then, since $\varphi _2$ is non-decreasing, we can find 
a $\delta > 0$ such that for all $u < u_0$ we have $\varphi _2 (u) \leq \varphi _2(u_0) -\delta.$
Also, by the definition of $\varphi _2$, there is $v>0$ such that 
$\varphi_2( u_0) \leq \varphi (u_0 v) -\varphi _1(v) +\frac{\delta }{3}$ and, by the left-continuity 
of $\varphi$, there is $t < u_0$ such that $0 \leq \varphi ( u_0 v) - \varphi (t v) \leq \frac{\delta }{3}$. 
Thus
\vspace{-3mm} 
\begin{equation*}
\varphi _2( t ) \geq \varphi ( t v) -\varphi_1( v) \geq \varphi ( u_0 v) - \varphi _1(v) -\frac{\delta }{3} 
\geq \varphi _2 ( u_0) -\frac{2\delta }{3},
\end{equation*}
which is a contradiction. This contradiction shows that $\varphi_2$ is left-continuous at $u_0 > 0$.
\smallskip

$2^{\small 0}.$ Let $ \varphi _2( u_0) =\infty $. Suppose again that $\varphi _2$ is 
not left-continuous at $u_0$. Then, since $\varphi _{2}$ is non-decreasing, we can find 
$M > 0$ such that for all $u < u_0$ we have that $\varphi _2( u) \leq M.$
Moreover, by the definition of $\varphi _2$, there is $v>0$ such that 
$\varphi ( u_0 v) - \varphi _1 (v) \geq 3M$ and, by the left-continuity of $\varphi $, there is  
$t < u_0$ such that $\varphi (t v) = \infty $ (in the case $u_0 v > b_{\varphi }$) or $\varphi ( t v)
\geq \varphi ( u_0 v) -M$ (in the case $u_0 v \leq b_{\varphi }$). Then, in the case 
$u_o v \leq b_{\varphi }$, we have
\begin{equation*}
\varphi _2( t) \geq \varphi ( t v) -\varphi_1( v) \geq \varphi ( u_0 v) -\varphi _1(v) -M \geq 2M,
\end{equation*}
or in the case $u_0 v > b_{\varphi }$ we obtain 
\begin{equation*}
\varphi _2 ( t ) \geq \varphi ( t v) -\varphi_1( v) = \infty \geq 3M,
\end{equation*}
which give contradictions. Thus, our claim is proved.
\medskip

(ii) By the definition of $\varphi_2$ we have $\varphi(uv) \leq \varphi_1(v) + \varphi_2(u)$ for 
all $u, v > 0$. Then (ii) follows from remarks after Corollary 2.
\medskip

(iii) The equivalence of the function $\frac{\varphi_1(u)}{\varphi(u v)}$ to a non-decreasing function 
means that there is a number $K > 0$ such that for each $v > 0$ 
there is a non-decreasing function $\psi_v$  with estimates 
$\frac{1}{K} \psi_v(u) \leq \frac{\varphi_1(u)}{\varphi(u v)} \leq K \psi_v(u)$ for all $u > 0$.

Let $u > 0$ be fixed and suppose $0 < \varphi_1^{-1}(u) < v$. Then, by the monotonicity of 
$\psi_w$, one has for $w = \frac{\varphi^{-1}(u)}{\varphi_1^{-1}(u)}$
\vspace{-2mm}
\begin{eqnarray*}
\frac{\varphi_1(v)}{\varphi(v w)}  
&\geq&
\frac{1}{K}\, \psi_w(v) \geq \frac{1}{K} \, \psi_w (\varphi_1^{-1}(u)) \geq
\frac{1}{K^2} \frac{\varphi_1(\varphi_1^{-1}(u))}{\varphi (\varphi_1^{-1}(u) w)} \\
&=& 
\frac{1}{K^2} \frac{u}{\varphi (\varphi^{-1}(u))} = \frac{1}{K^2},
\end{eqnarray*}
which gives $\varphi (v w) \leq K^2 \varphi_1(v)$. If $\varphi_1^{-1}(u) \geq v$, then by monotonicity 
of $\varphi$, we obtain $\varphi (v w) \leq \varphi (\varphi_1^{-1}(u) w) = \varphi (\varphi^{-1}(u)) = u$. 
Consequently, by convexity of $\varphi_1$, for any $v > 0$ we have that 
$\varphi (vw) \leq K^2 \varphi_1(v) + u \leq \varphi_1(K^2 v) + u$ and, therefore, 
$\varphi_2(\frac{w}{K^2}) \leq u$. Thus $\varphi^{-1}(u) \leq K^2 \varphi_1^{-1}(u) \varphi_2^{-1}(u)$ 
for all $u > 0$ and the proof is complete. \qed

\vspace{3mm} 
Using Lemma 5(iii), we obtain the following other version of Theorem 6.

\vspace{3mm} 
{\bf THEOREM 7}. \label{thm6} {\it Let $E$ be a Banach function
space with the Fatou property and let $\varphi, \varphi _{1}, \varphi _{2}$ be
Orlicz functions. Suppose 
\vspace{-2mm}
\begin{equation*}
M( E_{\varphi _1}, E_{\varphi }) \hookrightarrow E_{\varphi_2}.
\end{equation*}
Assume that for any $v>0$ the function $\frac{\varphi _1( u) }{\varphi (u v) 
}$ is equivalent to a non-decreasing function of $u > 0$. If $L^{\infty}
\hookrightarrow E$ and $E_{a} \not= \{ 0\} $, then $\varphi ^{-1}\prec
\varphi _{1}^{-1}\varphi _{2}^{-1}$ for large arguments.}

\vspace{3mm} 
{\it Proof.} By Lemma 5(iii), we know that there is $\varphi _{3}=\varphi
\ominus \varphi _{1}$ satisfying $\varphi _{1}^{-1}\varphi _{3}^{-1}\approx
\varphi ^{-1}$. Therefore, according to Corollary 4, we have 
\begin{equation*}
M(E_{\varphi _{1}}, E_{\varphi })=E_{\varphi _{3}}\hookrightarrow E_{\varphi
_{2}}.
\end{equation*}
Moreover, it is known (see \cite{KMP03}, Theorem 2.4) that if 
$E_{a}\not=\{0\} $ and $E_{\varphi _{3}}\hookrightarrow E_{\varphi _{2}}$
then there is $k>0$ such that $\limsup_{u\rightarrow \infty }\frac{\varphi
_{2}(k\,u)}{\varphi _{3}(u)}<\infty $. Therefore, we have $\varphi
_{2}(k\,u)\leq C\,\varphi _{3}(u)$ for some $C>1$ and large u. Consequently,
for $u=\varphi _{3}^{-1}(v)$ from Lemma 2 we obtain
\begin{equation*}
\varphi _{2}(k\varphi _{3}^{-1}(v))\leq C\,\varphi _{3}(\varphi
_{3}^{-1}(v))\leq Cv
\end{equation*}
and 
\begin{equation*}
k\varphi _{3}^{-1}(v)\leq \varphi _{2}^{-1}(\varphi _{2}(k\varphi
_{3}^{-1}(v)))\leq \varphi _{2}^{-1}(Cv)\leq C\varphi _{2}^{-1}(v) ~~ {\rm for ~large} ~v.
\end{equation*}
Finally, we have $\varphi ^{-1}\approx \varphi _{1}^{-1}\varphi
_{3}^{-1}\prec \varphi _{1}^{-1}\varphi _{2}^{-1}$ for large arguments and
the theorem is proved. \qed

\medskip 
It is worth to notice that there are Orlicz spaces $L^{\varphi},
L^{\varphi _1}$ such that for any $v>0$ the function $\frac{\varphi _1(u)}{
\varphi ( uv) }$ is non-decreasing in u, but there is no $a>0$ such that $
\frac{f_{L^{\varphi }}( t) }{f_{L^{\varphi_1}}( t)\, t^{a}}$ is
non-decreasing in t.

\vspace{3mm} 
{\bf Example 8}. Consider the Orlicz functions $\varphi (u) = u^2$ and $ \varphi
_1(u) = u^2 \ln(u+1)$. Then $L^{\varphi _{1}}[ 0, 1] \hookrightarrow L^{\varphi }[ 0, 1] $ 
and the function $\frac{\varphi_1( u) }{\varphi ( u v) }= \frac{\ln ( u+1) }{v^{2}}$ is 
non-decreasing in $u > 0$ for any $v>0$. On the other hand, if in the quotient 
$\frac{f_{L^{\varphi }}( t) }{f_{L^{\varphi _1}}( t) \, t^a} = \frac{\varphi
_{1}^{-1}( \frac{1}{t}) }{\varphi ^{-1}( \frac{1}{t}) \, t^a}$ after substitution 
$t = \frac{1}{\varphi _1( u) }$ we obtain 
\vspace{-2mm}
\begin{equation*}
\frac{\varphi _1( u)^a \, u}{\varphi^{-1}( \varphi_{1}( u)) } = \frac{
u^{2a+1}\ln ^a ( u+1) }{\sqrt{u^{2} \ln ( u+1) }} = u^{2a} \ln^{a-1/2} (
u+1) \rightarrow \infty ,
\end{equation*}
as $u \rightarrow \infty $ for any $a>0$. Consequently, $\frac{
f_{L^{\varphi}} ( t) }{f_{L^{\varphi _1}} ( t) \,t^a}\rightarrow \infty $ as 
$t\rightarrow 0^{+}$ and therefore it cannot be non-decreasing for small $t
> 0$.

\medskip 
If we drop the assumption that $\frac{\varphi _{1}(u)}{\varphi (uv)}
$ is non-decreasing in Theorem 7, then the result may not be true.

\vspace{3mm} 
{\bf Example 9}. Let $\varphi ( u) =\frac{u^2}{2}$ and we will construct
a new function $\psi$ which does not satisfy the $\Delta _2$-condition for
large arguments, i.e, $\lim \sup_{u \rightarrow \infty} \frac{\psi(2u)}{\psi(u)} = \infty$ 
and such that 
\begin{equation*}
\psi (u) \geq \varphi (u) ~ \mathrm{for ~all} ~~ u > 0 ~ \mathrm{and} ~~
\psi ( u_n) = \varphi (u_n)
\end{equation*}
for some sequence $( u_n) $ tending to infinity with $\frac{\psi ( 2 u_n) }{
\psi ( u_n) }\nearrow \infty$.

Take any sequence $( a_n)$ of positive real numbers satisfying two conditions 
\vspace{-2mm}
\begin{equation}  \label{sequence20}
\frac{a_{n+1}}{a_n}\nearrow \infty ~~ \mathrm{and} ~~ 2 \, \sum_{k=1}^{n} (
-1)^{n-k}\, a_k <\sum_{k=1}^{n+1}( -1)^{n+1-k} \,a_k ~~ \mathrm{for ~ all}
~~ n \in \mathbb{N}.
\end{equation}
It is easy to see that, for example, the sequence $a_{n}=(n+2)!$ satisfies those
conditions. Define the required sequence as $u_n = 2\, \sum_{k=1}^{n}
(-1)^{n-k} \, a_{k}, u_0 = 0$ and consider the sequence of pairwise disjoint
subintervals of $[0, \infty)$ defined by $I_n = [u_{n-1}, u_n), n = 1, 2,
\ldots$. The numbers $a_n$ are the centers of $I_n$, since $\frac{u_n + u_{n-1}}{2
} = a_n$. Now define the following Orlicz function  
\begin{equation}  \label{function21}
\psi ( u) = \int_{0}^{u} \sum_{n=1}^{\infty} a_{n} \,\chi _{I_n} (s) ds.
\end{equation}
For any $n \in \mathbb{N}$ we have 
\begin{equation*}
\int_{I_n} a_n ds = a_n(u_n - u_{n-1}) = \frac{1}{2} (u_n + u_{n-1})(u_n -
u_{n-1}) = \frac{u_n^2 - u_{n-1}^2}{2} = \int_{I_n} s ds
\end{equation*}
and, thus, 
\vspace{-2mm}
\begin{eqnarray*}
\psi(u_n) &=& \int_0^{u_n} \sum_{k=1}^{\infty} a_{k} \,\chi _{I_k} (s) ds =
\sum_{k=1}^{n} \int_{I_k} a_{k} ds \\
&=& \sum_{k=1}^{n} \int_{I_k} s ds = \int_0^{u_n} s ds = \frac{u_n^2}{2} =
\varphi(u_n).
\end{eqnarray*}
We must now check that the function $\psi $ is bigger than the function $\varphi .$ 
For $u\in \lbrack 0,u_{1}]=[0,2a_{1}]$ we have $\psi (u)=a_{1}u\geq u^{2}/2 $ and, 
for $u\in \lbrack u_{n-1},u_{n}],n=2,3,\ldots $, it yields 
\vspace{-2mm}
\begin{eqnarray*}
\psi (u) &=&\int_{0}^{u}\sum_{k=1}^{\infty }a_{k}\,\chi
_{I_{k}}(s)ds=\sum_{k=1}^{n-1}a_{k}(u_{k}-u_{k-1})+a_{n}\,(u-u_{n-1}) \\
&=&\frac{1}{2}\sum_{k=1}^{n-1}(u_{k}^{2}-u_{k-1}^{2})+\frac{u_{n}+u_{n-1}}{2}
(u-u_{n-1}) \\
&=&\frac{1}{2}\,u_{n-1}^{2}+\frac{u_{n}+u_{n-1}}{2}\,u-\frac{u_{n}+u_{n-1}}{2
}u_{n-1} \\
&=&\frac{u_{n}+u_{n-1}}{2}u-\frac{u_{n}u_{n-1}}{2}=\frac{h(u)}{2}+\frac{u^{2}
}{2},
\end{eqnarray*}
where 
\begin{equation*}
h(u)=-u^{2}+(u_{n}+u_{n-1})\,u-u_{n}\,u_{n-1}.
\end{equation*}
Since for $u\in \lbrack u_{n-1},u_{n}]$ one has
$h(u)\geq \max \, [h(u_{n-1}),h(u_{n})]=0$ it follows that $\psi (u)\geq \frac{u^{2}}{2}$ 
for any $u\in \lbrack u_{n-1},u_{n}]$, and consequently $\psi (u)\geq \frac{u^{2}
}{2}$ for any $u\geq 0$. Moreover, by assumptions (\ref{sequence20}) on 
$a_{n}$, we see that $2u_{n}\in I_{n+1}=[u_{n},u_{n+1})$ and one has 
\begin{eqnarray*}
\frac{\psi (2u_{n})}{\psi (u_{n})} &=&\frac{(u_{n+1}+u_{n})\,u_{n}-\frac{
u_{n+1}\,u_{n}}{2}}{u_{n}^{2}/2} \\
&=&\frac{2u_{n+1}\,u_{n}+2u_{n}^{2}-u_{n+1}u_{n}}{u_{n}^{2}}=2+\frac{u_{n+1}
}{u_{n}} \\
&=&2+\frac{2a_{n+1}-u_{n}}{u_{n}}=1+\frac{2a_{n+1}}{u_{n}} \\
&=&1+\frac{2a_{n+1}}{2a_{n}-u_{n-1}}>1+\frac{a_{n+1}}{a_{n}}\rightarrow
\infty ,
\end{eqnarray*}

\vspace{-3mm}

\noindent
as $n\rightarrow \infty $.

Of course, $L^{\psi }[0,1]\subset L^{\varphi }[0,1]=L^{2}[0,1]$ because 
$\psi (u)\geq \varphi (u)$ for all $u>0$ and thus $M(L^{\psi }, L^{\varphi })$
is non-trivial. Moreover, $L^{\psi }\not=L^{\varphi }=L^{2}$ since $\psi $
does not satisfy the $\Delta _{2}$-condition for large u. Let us calculate 
$\varphi _{2}=\varphi \ominus \psi $. For $u>1$ 
\begin{eqnarray*}
\varphi _{2}(u) &=&\sup_{v>0}\left[ \varphi (uv)-\psi (v)\right] \geq
\limsup_{n\rightarrow \infty }\left[ \varphi (uu_{n})-\psi (u_{n})\right]  \\
&=&\limsup_{n\rightarrow \infty }\frac{1}{2}u_{n}^{2}(u^{2}-1)=\infty ,
\end{eqnarray*}
and for $0<u\leq 1$ one has $\varphi (u\,v)-\psi (v)\leq \varphi (v)-\psi
(v)\leq 0$ for each $v>0$ and so $\varphi _{2}(u)=0$. 
Therefore, 
\begin{equation*}
~~\varphi _{2}(u)=
\begin{cases}
0 & ~\mathrm{if}~0\leq u\leq 1, \\ 
\infty  & ~\mathrm{if}~u>1,
\end{cases}
\end{equation*}
and $\psi ^{-1}(u)\varphi _{2}^{-1}(u)=\psi ^{-1}(u)\leq \varphi ^{-1}(u)$
for all $u>0$. 

\vspace{3mm}
\noindent
Let us now collect some properties of functions and spaces that arise in Example 9.

\begin{itemize}
\item[$(a)$] We don't have the relation $\varphi ^{-1}\prec \psi
^{-1}\varphi _{2}^{-1}$ for large u, since
$$
\liminf_{u\rightarrow \infty }\frac{\psi ^{-1}(u)}{\varphi ^{-1}(u)}
= \liminf_{v\rightarrow \infty }\frac{v}{\varphi ^{-1}(\psi (v))} 
$$
$$
\leq \lim_{n\rightarrow \infty }\frac{2u_{n}}{\varphi ^{-1}(\psi (2u_{n}))}
=\lim_{n\rightarrow \infty }\frac{\sqrt{2}u_{n}}{\sqrt{\psi (2u_{n})}} = 
\sqrt{2}\lim_{n\rightarrow \infty }\frac{u_{n}}{\sqrt{
(u_{n+1}+u_{n})u_{n}-u_{n+1}u_{n}/2}} 
$$
$$
=\sqrt{2}\lim_{n\rightarrow \infty }\frac{u_{n}}{\sqrt{
2a_{n+1}u_{n}-(2a_{n+1}-u_{n})u_{n}/2}} 
= \sqrt{2}\lim_{n\rightarrow \infty }\frac{u_{n}}{\sqrt{
a_{n+1}u_{n}+u_{n}^{2}/2}}
$$
$$
=\sqrt{2}\lim_{n\rightarrow \infty }\frac{1}{\sqrt{
\frac{1}{2}+\frac{a_{n+1}}{u_{n}}}} 
\leq \sqrt{2}\lim_{n\rightarrow \infty }\frac{1}{\sqrt{\frac{1}{2}+\frac{
a_{n+1}}{2a_{n}}}}=0.
$$

\item[$(b)$]  The function $\frac{\psi (u)}{\varphi (u)}$ cannot be monotone
since 
\begin{equation*}
\frac{\psi (u_{n})}{\varphi (u_{n})}=1~\mathrm{and}~~\frac{\psi (2u_{n})}{
\varphi (2u_{n})}=\frac{\psi (2u_{n})}{2u_{n}^{2}}=\frac{\psi (2u_{n})}{
4\psi (u_{n})}\rightarrow \infty ~~\mathrm{as}~n\rightarrow \infty .
\end{equation*}
\end{itemize}

\noindent
In consequence, from (a) and (b), we obtain 
\begin{itemize}
\item[$(c)$] we cannot drop the assumption of monotonicity of $\frac{\psi
(u)}{\varphi (uv)}$ in Lemma 5(iii), in general.
\end{itemize}

\noindent
Moreover,
\begin{itemize}
\item[$(d)$] there is no $a>0$ such that $\frac{f_{L^{\varphi }}(t)}{f_{L^{\psi
}}(t)t^{a}}$ is non-decreasing near zero, because
\begin{equation}
\limsup_{t\rightarrow 0^{+}}\frac{f_{L^{\varphi }}(t)}{f_{L^{\psi }}(t)}
=\limsup_{u\rightarrow \infty }\frac{\psi ^{-1}(u)}{\varphi ^{-1}(u)}\geq
\limsup_{n\rightarrow \infty }\frac{\psi ^{-1}(u_{n})}{\varphi ^{-1}(u_{n})}
=1  \label{equal29}
\end{equation}
and for each $a>0$ and every sequence $t_{n}\rightarrow 0^{+}$ we have 
$t_{n}^{-a}\rightarrow +\infty .$

\item[$(e)$] The function $\frac{f_{L^{\varphi }}(t)}{f_{L^{\psi }}(t)}$ cannot be
equivalent at $0$ to any pseudo-concave function, because of (\ref{equal29})
and 
\begin{equation*}
\liminf_{t\rightarrow 0^{+}}\frac{f_{L^{\varphi }}(t)}{f_{L^{\psi }}(t)}
=\liminf_{t\rightarrow 0^{+}}\frac{\psi ^{-1}(\frac{1}{t})}{\varphi ^{-1}(
\frac{1}{t})}=\liminf_{u\rightarrow \infty }\frac{\psi ^{-1}(u)}{\varphi
^{-1}(u)}=0.
\end{equation*}
Thus, in particular, formula (5.21) $f_{M(E,F)}(t)=\frac{f_{F}(t)}{f_{E}(t)}$ in the book 
\cite{AZ90} is false, in general (even up to equivalence). Note that in this example 
we have $f_{M(E, F)}(t) = \sup_{0 < s \leq t} \frac{f_F(s)}{f_E(s)} = 1.$

\item[$(f)$] We have $M(L^{\psi }, L^{\varphi }) = L^{\infty }=L^{\varphi _{2}}.$ 

In fact, as we have seen in the proof of Theorem 1(iii), we always have that 
$f_{M(E, F)}(t)\geq \frac{f_{F}(t)}{f_{E}(t)}$. Therefore, 
$ \limsup_{t\rightarrow 0^{+}}f_{M(L^{\psi }, L^{\varphi })}(t)\geq
\limsup_{t\rightarrow 0^{+}}\frac{f_{L^{\varphi }}(t)}{f_{L^{\psi }}(t)}\geq
1$ and $M(L^{\psi }, L^{\varphi })$ as a symmetric space on $[0,1]$ with the Fatou
property such that 

$\lim_{t\rightarrow 0^{+}}f_{M(L^{\psi }, L^{\varphi
})}(t)>0$ must be $L^{\infty }[0,1]$ (cf. \cite{LT79}, p. 118).

\item[$(g)$] The assumption of Theorem 6(i) is not satisfied since we have (d) and 
also the assumption of Theorem 7 is not satisfied since we have (b) but one can
see that $M(L^{\psi }, L^{\varphi })\subset L^{\varphi _{2}}$ while $\varphi
^{-1}\not\prec \psi ^{-1}\varphi _{2}^{-1}$ as it was shown in (a).

\item[$(h)$] There is no Young function $\varphi_3$ satisfying the equivalence 
$\varphi^{-1} \approx \psi^{-1} \varphi_3^{-1}$. If such a function should exists, then, 
by Corollary 4, we have that $M(L^{\psi}, L^{\varphi}) = L^{\varphi_3}$. On the other hand, from 
(f) we have $M(L^{\psi}, L^{\varphi}) = L^{\infty}$, which will mean that 
$\varphi_3 \approx \varphi_2$ for large arguments, but this is not possible because of (a).

\end{itemize}

\noindent
As we have seen earlier we have equality $M(E, E) = L^\infty$ and
we can ask if $M(E, F) = L^\infty$ implies that $E = F$?  From (f) we see that this is not always the case.
\vspace{1mm}

Assume here that ${\it supp} E = {\it supp} F = \Omega$. Note that $M(E, F) = L^\infty$ if and only if 
$E^{F F} = F$. Really, if $E^F = L^\infty$, then $E^{F F} = (L^{\infty})^F = F$. On the other hand, if 
$E^{F F} = F$, then $E^F = E^{F F F} = F^F = L^{\infty}$. We can also have
equality $M(E, F) = L^\infty$ if $E$ is a proper subspace of $F$ and the
norms of $\|\cdot\|_E$ and $ \| \cdot\|_F$ are equivalent on $E$. 
For example, if $E = F_a$ with $F \neq F_a$ and ${\it supp} F_a = {\it supp} F = \Omega$. 
Thus
\begin{itemize}
\item[$(i)$]  {\it $L^{\psi }$ is not $L^{\varphi }$-perfect}. In fact, from (f) and 
$L^{\psi } \neq L^{\varphi }$ we obtain 
\begin{equation*}
(L^{\psi })^{L^{\varphi} L^{\varphi }} = M(M(L^{\psi }, L^{\varphi }),
L^{\varphi }) = M(L^{\infty}, L^{\varphi }) = L^{\varphi } \neq L^{\psi }.
\end{equation*}
\end{itemize}

Using Lemma 5, Corollary 4(ii) and the operation $\varphi \ominus \varphi _1$ we are able to prove that 
the multiplier space between two Orlicz spaces $M(L^{\varphi_1}, L^{\varphi})$ on $[0, 1]$ is an Orlicz space 
$L^{\varphi_2}$ with $\varphi_2 = {\varphi \ominus \varphi _1}$ under some additional assumptions on 
the Orlicz functions $\varphi, \varphi_1$, which is a certain similarity to the case of sequence Orlicz spaces. 
Our proof is presented even for the Calder\'{o}n-Lozanovski\u{\i} spaces.

\newpage
{\bf Theorem 8}. \label{cor-6} {\it Let $\varphi, \varphi_1$ be increasing Orlicz functions and let $E$ 
be a symmetric space on $[0, 1]$ with the Fatou property.
\vspace{-3mm} 

\begin{itemize}
\item[$(i)$] If $\limsup_{u \rightarrow \infty} \frac{\varphi(u v)}{\varphi_1(u)} = 0$ for any $v > 0$ and additionally 
at least one of the following three conditions holds: either the function $f_v(u): = \frac{\varphi(uv)}{\varphi_1(u)}$ is 
non-increasing on $(0, \infty)$ for any $v > 0$ or $\frac{\varphi^{-1}(u)}{\varphi_1^{-1}(u)}$ is a non-decreasing 
function for large $u$ or the function 
\begin{equation} \label{Thm8equation28}
\varphi_2 (u) = (\varphi \ominus \varphi _{1})(u) = \sup_{v > 0}
[\varphi(uv) - \varphi_1(v)]
\end{equation}

\vspace{-4mm}

satisfies the $\Delta_2$-condition for large arguments, then $M(E_{\varphi_1}, E_{\varphi}) = E_{\varphi_2}$.

\item[$(ii)$] If $\limsup_{u \rightarrow \infty} \frac{\varphi(u v)}{\varphi_1(u)} <
\infty $ for some $v > 0$ and $\limsup_{u \rightarrow \infty} \frac{
\varphi(u w)}{\varphi_1(u)} > 0 $ for some $w > 0$, then
$M(E_{\varphi_1}, E_{\varphi}) = L^{\infty}$.
\vspace{-2mm}

\item[$(iii)$] If $\limsup_{u \rightarrow \infty} \frac{\varphi(u v)}{\varphi_1(u)} =
\infty$ for all $v > 0$, then $M(E_{\varphi_1}, E_{\varphi}) = \{ 0\}$.
\end{itemize}
}

{\it Proof.} (i) We only need to prove that in all these three cases we have $\varphi
^{-1}\prec \varphi _{1}^{-1}\varphi _{2}^{-1}$ for large arguments, since 
by Lemma 5(ii), we have $\varphi _{1}^{-1}\varphi _{2}^{-1} \prec \varphi^{-1}$ 
even for all arguments, which means that 
$\varphi _{1}^{-1}\varphi _{2}^{-1}\approx \varphi ^{-1}$ for large arguments and 
then Corollary 4(ii)l implies $M(E_{\varphi_1}, E_{\varphi}) = E_{\varphi_2}$.
Therefore, in each of these three cases we will proceed as follows:

$1^{\small 0}$. If for any $v > 0, f_v(u)$ is a non-increasing function on $(0, \infty)$, 
then, by Lemma 5(iii), we obtain that $\varphi _{1}^{-1}\varphi _{2}^{-1}\approx \varphi ^{-1}$ 
for all arguments. 
\smallskip

$2^{\small 0}$. Let $\frac{\varphi^{-1}(u)}{\varphi_1^{-1}(u)}$ be non-decreasing for $u > u_0 \geq 0$. 
Since $\limsup_{v \rightarrow \infty} \frac{\varphi(v w)}{\varphi_1(v)} = 0$ for any $w > 0$ it follows 
that the supremum in the definition of $\varphi_2$ is attained at some $v_0 = v_0(w) > 0$. 
Can we say something more about this $v_0$? If $w = \frac{\varphi^{-1}(u)}{\varphi_1^{-1}(u)}$ 
and we have $v_0 \geq \varphi_1^{-1}(u)$, then, by the monotonicity assumption, we get 
$$
\varphi_2[\frac{\varphi^{-1}(u)}{\varphi_1^{-1}(u)}] = \varphi [\frac{\varphi^{-1}(u)}{\varphi_1^{-1}(u)}\, v_0] 
- \varphi_1(v_0) \leq \varphi [\frac{\varphi^{-1}(\varphi_1(v_0))}{v_0}\, v_0] - \varphi_1(v_0) = 0
$$ 
and this case is not important since $\varphi_2 \geq 0$. Therefore it must be $v_0 \leq \varphi_1^{-1}(u)$ 
with $u > u_0 \geq 0$, which, in its turn gives,
$$ 
\varphi_2[\frac{\varphi^{-1}(u)}{\varphi_1^{-1}(u)}] = \varphi [\frac{\varphi^{-1}(u)}{\varphi_1^{-1}(u)}\, v_0] 
- \varphi_1(v_0) \leq \varphi [\frac{\varphi^{-1}(u)}{\varphi_1^{-1}(u)}\, v_0] \leq \varphi [\varphi^{-1}(u)] = u,
$$
i.e., $\varphi^{-1}(u) \leq \varphi _{1}^{-1}(u)\, \varphi _{2}^{-1}(u)$ for $u > u_0 \geq 0$.
\smallskip

$3^{\small 0}$. Let $\varphi_2$ satisfy the $\Delta_2$-condition for large arguments, that is, there exist 
constants $C \geq 1, u_0 \geq 0$ such that $\varphi_2(2 u) \leq C \varphi_2(u)$ for all $u > u_0$. 
Similarly as in $2^{\small 0}$ we find that for any $w > 0$ there exists a $v_0 = v_0(w) > 0$ such 
that $\varphi_2(w) = \varphi(w v_0) - \varphi_1(v_0)$. For $w = \varphi_2^{-1}(u)$ we have 
$$
u = \varphi_2[\varphi_2^{-1}(u)] = \varphi [ \varphi_2^{-1}(u)\, v_0] - \varphi_1(v_0) > 0,
$$
that is, $\varphi_2^{-1}(u) \geq \frac{\varphi^{-1}[\varphi_1(v_0)]}{v_0}$. Hence, by using Lemma 5(ii), we obtain 
$$
1 \geq \frac{\varphi^{-1}[\varphi_1(v_0)]}{v_0 \, \varphi_2^{-1}(u)} \geq 
\frac{\varphi_1^{-1}[\varphi_1(v_0)]\, \varphi_2^{-1}[\varphi_1[v_0)]}{2 v_0\, \varphi_2^{-1}(u)} 
= \frac{\varphi_2^{-1}[\varphi_1(v_0)]}{2\, \varphi_2^{-1}(u)},
$$
and, by the $\Delta_2$-condition of $\varphi_2$ for $u > u_1 = \varphi_2(u_0)  \geq 0$, we get 
$$
v_0 \leq \varphi_1^{-1} [ \varphi_2(2 \, \varphi_2^{-1}(u))] \leq \varphi_1^{-1}(C u).
$$
Since $\varphi^{-1}[u + \varphi_1(v_0)] = \varphi_2^{-1}(u)\, v_0$ it follows that
\begin{eqnarray*}
\varphi^{-1}(u) 
&\leq&
\frac{\varphi_1^{-1}(Cu)}{v_0} \varphi^{-1}(u) \leq \frac{\varphi_1^{-1}(Cu)}{v_0} \varphi^{-1}[u + \varphi_1(v_0)]\\
&\leq&
\frac{\varphi_1^{-1}(Cu)}{v_0} \varphi_2^{-1}(u) v_0 = \varphi_1^{-1}(C u)\, \varphi_2^{-1}(u) \leq 
C \varphi_1^{-1}(u)\, \varphi_2^{-1}(u)
\end{eqnarray*}
for $u > u_1 = \varphi_2(u_0) \geq 0$. Therefore, all three cases are proved.
\vspace{2mm}

(ii) Suppose $\limsup_{u\rightarrow \infty }\frac{\varphi (uv)}{
\varphi _{1}(u)}<\infty $ for some $v>0$. Then there is $K>0$ such that 
$\varphi (uv)\leq K\varphi _{1}(u)$ for large $u$ and, by \cite[Theorem 2.3]{KMP03},
we obtain $E_{\varphi _{1}}\hookrightarrow E_{\varphi }$. 
Thus $L^{\infty }\hookrightarrow M(E_{\varphi _{1}}, E_{\varphi })$. 
On the other hand, suppose on the contrary that $\limsup_{u\rightarrow \infty }\frac{\varphi (uw)}{\varphi
_{1}(u)}=\eta >0$ for some $w>0$ and $M(E_{\varphi _{1}}, E_{\varphi})\not=L^{\infty }$. 
Define the new function $\psi \left( u\right) =\frac{2}{\eta }\varphi (uw)$.
Then, again, by \cite[Theorem 2.3]{KMP03} we have that $E_{\varphi } = E_{\psi }$ and so 
$M(E_{\varphi _{1}}, E_{\varphi }) = M(E_{\varphi _{1}}, E_{\psi
})$. The fundamental function $f_{M}$ of the symmetric space $M= M(E_{\varphi _{1}}, E_{\psi
})$ satisfies the condition $\lim_{t\rightarrow 0^{+}} f_{M}( t) = 0$ because 
$M(E_{\varphi _{1}}, E_{\psi })\neq L^{\infty }$. Since
\begin{equation*}
1 = \left\| \frac{\chi _{[0, t]}} {f_{M}(t) }\right\|_{M(E_{\varphi _{1}}, E_{\psi })} 
\geq \left\| \frac{\chi _{[0, t]}} {f_{M}( t) } \frac{\chi_{[0, t]}} {f_{E_{\varphi _{1}}}( t) }
\right\|_{\psi } =
\frac{1}{f_{M}( t) }\frac{f_{E_{\psi }}( t) }{ f_{E_{\varphi _{1}}}( t) }
\end{equation*}
and $\lim_{t \rightarrow 0^+} f_{M}( t) = 0$ it follows that $\lim_{t \rightarrow 0^+} \frac{f_{E_{\psi
}} ( t) }{f_{E_{\varphi _{1}}}( t) } =  0$. This means 
\begin{equation*}
0=\lim_{t\rightarrow 0^{+}}\frac{f_{E_{\psi }}( t) }{
f_{E_{\varphi _{1}}}( t) }  = \lim_{t \rightarrow 0^+} \frac{\varphi _{1}^{-1}( 1/f_E(t)) }{\varphi^{-1}( 1/f_E(t)) } 
= \lim_{u\rightarrow \infty }\frac{\varphi _{1}^{-1}( u) }{\psi ^{-1} ( u) }.
\end{equation*}
But 
\begin{equation*}
\eta <\limsup_{u\rightarrow \infty }\frac{\varphi ( uw) }{\varphi
_{1}( u) }=\limsup_{u\rightarrow \infty }\frac{\frac{\eta }{2}
\psi( u) }{\varphi _{1}( u) },
\end{equation*}
and thus we can find a sequence $u_{n}\rightarrow \infty $ such that 
$\psi ( u_{n}) \geq \varphi _{1}( u_{n}) $. Putting $v_{n}=\psi(u_{n}) $ we see that 
$\frac{\varphi _{1}^{-1}\left( v_{n}\right) }{\psi ^{-1}\left( v_{n}\right) }\geq 1$, which is 
a contradiction with the just mentioned equality 
$\lim_{u\rightarrow \infty }\frac{\varphi _{1}^{-1}( u) }{\psi^{-1}( u) }=0$.
\vspace{0.1mm}

(iii) The condition: there are $K, u_0, M > 0$ such that $\varphi ( Ku) \leq
M \varphi _{1}( u)$ for all $u > u_0$ is necessary for the inclusion $E_{\varphi _{1}}\hookrightarrow E_{\varphi }$
(see \cite{KMP03}, Theorem 2.4) and this inclusion is necessary for $M(E_{\varphi
_{1}},E_{\varphi })\not=\{0\}$ by Proposition 1(i). But $\limsup_{u\rightarrow
\infty }\frac{\varphi (uv)}{\varphi _{1}(u)}=\infty $ means that the just mentioned condition 
on the function $\varphi$ is not satisfied. 

\vspace{3mm} 
{\bf Example 10}. Let $\varphi_1$ be an increasing Orlicz function and $\varphi(u) = 2 \varphi_1(\sqrt{u})$. 
Then $\varphi_2(u) = \sup_{v > 0} [ \varphi(uv) - \varphi_1(v)] = \varphi_1(u)$ since, by convexity of $\varphi_1$, 
we have 
$$
\varphi(uv) - \varphi_1(v) = 2 \varphi_1(\sqrt{uv}) - \varphi_1(v) \leq 2 \varphi_1(\frac{u+v}{2}) - \varphi_1(v) \leq \varphi_1(u)
$$
with equality for $v = u$ (see also \cite[p. 269]{ZR67} and \cite[p. 79]{Ma89}). If $\varphi$ is a convex function, then 
$$
\varphi^{-1}(u) = \varphi_1^{-1}(u/2)^2 \leq \varphi_1^{-1}(u)^2 = \varphi^{-1}(2u) \leq 2 \varphi^{-1}(u)
$$
for any $u > 0$ and from Corollary 4 we obtain that $M(E_{\varphi_1}, E_{\varphi}) = E_{\varphi_1}$ for any Banach 
ideal space $E$ with the Fatou property. Note that for the concrete $\varphi_1(u) = \exp(u^2) - 1$ we have that 
$\varphi(u) = 2(e^u-1)$ and $f_4(u) = \frac{\varphi(4u)}{\varphi_1(u)}$ is not decreasing on $(0, \infty)$ since 
$f_4^{\prime}(1) > 0$, $\varphi_1$ does not satisfy the $\Delta_2$-condition for large arguments, but the 
function $\varphi^{-1}(u)/\varphi_1^{-1}(u)$ is increasing on $(0, \infty)$.

\vspace{3mm}
Theorem 8 with $E = L^1[0, 1]$, that is, for Orlicz spaces $L^{\varphi}, L^{\varphi_1}$ on $[0, 1]$ and the space 
of multipiers $M(L^{\varphi_1}, L^{\varphi})$ has the following form:

\vspace{3mm} 
{\bf Corollary 6}. \label{cor-6} {\it Let $\varphi, \varphi_1$ be increasing Orlicz functions generating the corresponding 
Orlicz spaces $L^{\varphi}$ and $L^{\varphi_1}$ on $[0, 1]$.
\vspace{-3mm} 

\begin{itemize}
\item[$(i)$] If $\limsup_{u \rightarrow \infty} \frac{\varphi(u v)}{\varphi_1(u)} = 0$
for any $v > 0$ and additionally at least one of three conditions on $\varphi, \varphi_1, \varphi_2$ 
from Theorem 8(i) hold, then $M(L^{\varphi_1}, L^{\varphi}) = L^{\varphi_2}$.
\vspace{-2mm} 

\item[$(ii)$] If $\lim_{u \rightarrow \infty} \frac{\varphi(u v)}{\varphi_1(u)} <
\infty $ for some $v > 0$ and $\limsup_{u \rightarrow \infty} \frac{
\varphi(u w)}{\varphi_1(u)} > 0 $ for some $w > 0$, then
$M(L^{\varphi_1}, L^{\varphi}) = L^{\infty}$.
\vspace{-2mm}

\item[$(iii)$] If $\limsup_{u \rightarrow \infty} \frac{\varphi(u v)}{\varphi_1(u)} =
\infty$ for all $v > 0$, then $M(L^{\varphi_1}, L^{\varphi}) = \{ 0\}$.

\end{itemize}
}

Corollary 6(i) without any proof was written in the papers by Wang \cite[Lemma 2]{Wa63}, 
Zabre{\u \i}ko \cite[p. 109]{Za66} and in the book by Appell and Zabrejko \cite[p. 123]{AZ90}. 
In these mentioned sources the authors formulate Corollary 6(i) without additional assumptions 
on $\varphi, \varphi_1$, but we were able to prove only the result with these three additional 
assumptions. Of course, it will be nice to give the proof without these additional conditions. 
The proof of the first case in Corollary 6(i) was already given in the book \cite[pp. 77-78]{Ma89}. 
Note that Ando \cite[Theorem 5]{An60} for given Orlicz functions $\varphi_1, \varphi$ defined 
the function $\varphi_2$ by the formula (\ref{Thm8equation28}) and proved that $L^{\varphi_2}$ 
is a largest Orlicz space on $[0, 1]$ such that 
$L^{\varphi_2} \subset M(L^{\varphi_1}, L^{\varphi}) \neq \{0\}$. 

Parts (ii) and (iii) appeared without any proof in Zabre{\u \i}ko \cite[pp. 108-109]{Za66} and with different proofs than 
our in the book \cite[pp. 132, 148-149]{AZ90} (the same proof of part (iii) appeared also earlier in \cite[pp. 309-310]{AZ87}).

Already in 1957 Shragin [Sh57] proved that $x \in M(L^{\varphi_1}, L^{\varphi})$ if and only if there are $c > 0$ and 
$\lambda > 0$ such that $\int_0^1 \varphi( \lambda |x(t)| \, f_c(\lambda |x(t)|) \, dt < \infty$, where 
$f_c(u) = \sup\{ v \geq 0: \varphi(u v) \geq c \varphi_1(v)\}$. It seems that the last condition, in general, cannot be 
discribed in terms of the function $\varphi_2$.

B. Maurey in the paper [29] on pages 128-138 is proving that if $ \varphi, \varphi_1$ are two Orlicz functions 
which additionally are N-functions at infinity, that is, 
$\lim_{u \rightarrow \infty} \frac{\varphi(u)}{u} = \lim_{u \rightarrow \infty} \frac{\varphi_1(u)}{u} \newline = \infty$ 
and $\lim_{u\rightarrow \infty }\frac{\varphi ( vu) }{\varphi_1 ( u) }=0$ for any $v>0$, then for 
any measure space $(\Omega, \mu)$ we have (cf. [29], Proposition 107) 
\begin{equation*}
M \left ( L^{\varphi_1 }(\Omega, \mu), L^{\varphi }(\Omega, \mu)\right) = L^{\theta }(\Omega, \mu),
\end{equation*}
where $\theta = \varphi \ominus \varphi_1$.
His proof of this result is using one important property of operation $ \ominus$, namely that
$\varphi \ominus [\varphi \ominus \varphi_1] (u) = \varphi_1(u)$ for all $u > 0$ (cf. [29], Proposition 104(b), p. 130). 
Unfortunately, the last equality is not true for all $u > 0$, as we can see on the example below.

\vspace{3mm}
\textbf{Example 11.}  Let $\varphi ( u) =u^2$ and let $\varphi _p$ for $1 \leq p \leq 4$ be defined by
\begin{equation*}
~~\varphi _p(u) = 
\begin{cases}
u^p & ~\mathrm{if} ~0 \leq u \leq 1, \\ 
u^4 & ~\mathrm{if} ~u \geq 1.
\end{cases}
\end{equation*}
If $1 \leq p \leq 2$, then
\begin{equation*}
~~\theta (u) = \varphi \ominus \varphi _p(u) = 
\begin{cases}
0 & ~\mathrm{if} ~0 \leq u \leq 1, \\ 
u^2 -1 & ~\mathrm{if} ~1 \leq u \leq \sqrt{2},\\
u^4/4 & ~\mathrm{if} ~ u \geq \sqrt{2},
\end{cases}
\end{equation*}
and
\begin{equation*}
~~  \varphi \ominus \theta (u) =  \varphi \ominus [\varphi \ominus \varphi_p](u) = \varphi_2(u) =
\begin{cases}
u^2 & ~\mathrm{if} ~0 \leq u \leq 1, \\ 
u^4 & ~\mathrm{if} ~ u \geq 1.
\end{cases}
\end{equation*}
Therefore, for $1 \leq p < 2,  \varphi \ominus [\varphi \ominus \varphi_p]$ is equal to $\varphi_p$ only 
on interval $[1, \infty)$ but not on the interval $(0, 1)$ where it is $\varphi_2$.
\vspace{3mm}

We finish our considerations with a conjecture motivated by the above Example 9, Theorem 8 and Example 10.

\vspace{3mm} 
{\bf Conjecture.} We have equality $M(E_{\varphi _{1}}, E_{\varphi
})=E_{\varphi \ominus \varphi _{1}}$ for any Banach ideal space $E$.


\vspace{3mm}

\noindent {\footnotesize Pawe\l\ Kolwicz and Karol Le\'{s}nik, Institute of
Mathematics of Electric Faculty\newline
Pozna\'n University of Technology, ul. Piotrowo 3a, 60-965 Pozna\'{n}, Poland
}\newline
\textit{E-mails:} ~\texttt{pawel.kolwicz@put.poznan.pl,
klesnik@vp.pl}\newline

\vspace{-1mm}

\noindent {\footnotesize Lech Maligranda, Department of Engineering Sciences and Mathematics\\
Lule\aa\ University of Technology, SE-971 87 Lule\aa , Sweden}\newline ~\textit{E-mail:} \texttt{lech.maligranda@ltu.se
}

\end{document}